\documentclass[a4paper,12pt]{article}

\usepackage[austrian,american]{babel}

\usepackage[intlimits,leqno]{amsmath}
\usepackage{amssymb}
\usepackage{amsthm}
\usepackage{graphicx,color}

\theoremstyle{plain}
\newtheorem{*theo}{*Theorem}
\newtheorem{theo}{Theorem}
\newtheorem{prop}{Proposition}

\newtheorem{lemm}{Lemma}

\newtheorem{rema}{Remark}
\newtheorem{exam}{Example}
\newtheorem{defi}{Definition}

\renewcommand{\a}{a}
\newcommand{\A}{\mathcal{A}}

\newcommand{\B}{\mathcal{B}}

\renewcommand{\i}{{\rm i}}
\newcommand\R{{\mathbb{R}}}
\newcommand\N{{\mathbb{N}}}
\newcommand\C{{\mathbb{C}}}

\newcommand{\F}{\mathcal{F}}

\newcommand{\x}{\mathbf{x}}
\renewcommand{\k}{\mathbf{k}}
\newcommand{\y}{\mathbf{y}}
\newcommand{\z}{\mathbf{z}}

\newcommand\vj{{\mathbf{j}}}

\newcommand\vc{{\mathbf{c}}}

\newcommand\vO{{\mathbf{0}}}

\usepackage{makeidx}
\usepackage{epsfig}

\usepackage{latexsym}

\usepackage{pstricks,pst-node,amsfonts}

\newcommand{\const}{\mathop{\hspace{0.01pt}{\rm const}}\nolimits}
\renewcommand{\d}{{\rm d}}

\renewcommand{\S}{\mathcal{S}}
\newcommand{\supp}{\mbox{supp}}

\begin{document}

\title{On the causality of real-valued semigroups and diffusion}

\author{ Richard Kowar\\
Department of Mathematics, University of Innbruck, \\
Technikerstrasse 21a/2, A-6020,Innsbruck, Austria
}

\maketitle

\begin{abstract}
In this paper we show that a process modeled 
by a strongly continuous real-valued semigroup (that has a space convolution 
operator as infinitesimal generator) cannot satisfy causality. 
We present and analyze a causal model of diffusion that 
satisfies the semigroup property at a discrete set of time points 
$M:=\{\tau_m\,|\,m\in\N_0\}$ and that is in contrast to the classical diffusion 
model not smooth.
More precisely, if $v$ denotes the concentration of a substance diffusing 
with constant speed, then $v$ is continuous but its time derivative is 
discontinuous at the discrete set $M$ of time points. 
It is this property of diffusion that forbids the classical limit procedure
that leads to the noncausal diffusion model in Stochastics. 
Furthermore, we show that diffusion with constant speed satisfies an inhomogeneous 
wave equation with a time dependent coefficient. 
\end{abstract}

\section{Introduction}

The standard model of diffusion and its variants have many applications, for example, in 
solid state physics (\cite{KitKro93}), environmental modeling (\cite{Harris79}), image processing 
(\cite{GuiMorRya04,ScGrGrHaLe09,Soi99,Wei98a}) and inverse problems 
(\cite{EngRun95,EngHanNeu96,Isa98,Kir96,ScGrGrHaLe09}). 
For these applications the standard diffusion model (\cite{Fet80})
\begin{equation}\label{stddiffeq}
\begin{aligned}
     &\frac{\partial v}{\partial t}  - \nabla\cdot (D_0\,\nabla v) = 0
               \qquad \mbox{ $t>0$}\qquad \mbox{ with } \qquad v|_{t=0}=u \,
\end{aligned}
\end{equation}
is usually modified by replacing the diffusion constant $D_0$ by a tensor that takes inhomogenities 
or anisotropy of the medium into account. Recently \emph{sub-diffusion} (very slow diffusion) 
and \emph{super-diffusion} (very fast diffusion) have been modeled by a completely different approach 
where the time derivative in~(\ref{stddiffeq}) is replaced by a fractional time derivative 
(see e.g. \cite{KiSrTr06}). 
Although, many of these models (or their discretizations) are successfully used in applications, 
they violate the principle of causality. By causality we mean, for example, that a characteristic 
feature of a process like an \emph{interface} or a \emph{front} must propagate with a finite speed. 
It is peculiar that causality is a fundamental concept in physics of wave propagation and the 
respective mathematical theory of hyperbolic equations, but it is ignored in physics of diffusion 
and the respective mathematical theory of parabolic equations and stochastic processes. 
Is it so difficult - or even impossible - to find a causal model of diffusion?  
It seems that processes like diffusion are more ``complicated'' than processes related to waves, 
but as a matter of fact many of the standard models of wave propagation in \emph{dissipative media} 
(cf.~\cite{KowSchBon10,Kow10,KowSch10,NacSmiWaa90}) lack causality, too.  
Causality demands a coupling of space and time and therefore forbids \emph{model equations} 
that are to ``simple". 
Hence, apart from its physical significance, it is interesting to study causality from the mathematical 
point of view. 
Even if causality is not relevant for many applications, which is not really  clear yet, there is no 
mathematical reason to ignore this issue.  

It is the concern of this paper to model and analyze a causal model of diffusion.  
The paper consists of essentially two parts.  
First we show that causal diffusion cannot satisfy the strongly continuous semigroup property and 
its respective evolution equation. More precisely, we show that a process decribed by the 
\emph{scalar-valued} evolution equation
\begin{equation}\label{evolequv}
\begin{aligned}
     &\frac{\partial v}{\partial t}  - \A\,v = 0
               \qquad \mbox{ $t>0$}\qquad \mbox{ with } \qquad v|_{t=0}=u \,,
\end{aligned}
\end{equation}
where the space convolution operator $\A$ is the infinitesimal generator of a strongly continuous 
semigroup (\cite{DauLio92_5,EVa99}), cannot satisfy causality. That is to say, there exists a 
fundamental solution 
$G:\R^N\times\R\to\R$ of~(\ref{evolequv}) that does not satisfy 
\begin{equation}\label{causaGcond}
\begin{aligned}
   \supp (G) \subseteq \{(\x,t)\in\R^N\times [0,\infty)\,|\,|\x|\leq c_0\,t\}\,
\end{aligned}
\end{equation}
for a positive constant $c_0$. This causality condition can be interpreted as follows: 
Consider the diffusion of an ink droplet of radius $R_0$ dropped in water at time $t=0$.
Then causality demands that the ink water interface propagates with a finite speed bounded 
by $c_0$. Even if there is no interface (visible), the concentration $v$ of ink must vanish 
outside a closed ball of radius $R(t)\leq R_0+c_0\,t$.  

Although diffusion equation~(\ref{stddiffeq}) does not satisfy~(\ref{causaGcond}), 
many scientists consider this equation as \emph{weak causal} (or
shortly as \emph{causal}), since there exists a fundamental solution $G$ that satisfies
\begin{equation}\label{defcaus01}
   G(\cdot,t) = 0 \qquad\mbox{for}\qquad t<0\,.
\end{equation}
However, this property has nothing to do with causality. As a consequence, the problem of 
causality has beed ignored.

The second and larger part of this paper presents and analyzes a causal model of diffusion. 
In case of constant speed $c_0$ of diffusion, the model reads as follows
\begin{equation}\label{defcausdiff10}
\begin{aligned}
   v(\x,t)
      = \frac{\int_{S_{R(t)}(\x)} v\left(\x',n(t)\,\tau\right)
      \d \sigma(\x') }{|S_{R(t)}(\vO)|}\,
\quad\mbox{with}\quad v|_{t=0}=u\,,
\end{aligned}
\end{equation}
where $\tau>0$ denotes a small time period, $n(t)\in\N_0$ is such that $t\in (n\,\tau,(n+1)\,\tau]$ and 
$R(t):=c_0\,(t-n(t)\,\tau)$. Moreover, $\d \sigma(\x')$ denotes the Lebesgue surface measure on 
$\R^N$ and $|S_{R}(\vO)|$ denotes the surface area of the sphere $S_R(\vO)$. 
This diffusion modell satisfies causality, the semigroup property on the discrete set 
of time points $M=\{m\,\tau\,|\,m\in \N_0\}$ and $t\mapsto\frac{\partial v}{\partial t}(\x,t)$ is 
discontinuous at $t\in M$ for each $\x\in\supp\,(v(\cdot,t))$. 
Thus model~(\ref{defcausdiff10}) constitutes a compromise between causality and the strongly 
continuous semigroup property. 
Physically, this model describes a transport process of mass during which the masses 
split up at time points $t\in M$ and then spread out in each direction without 
interfering.\footnote{The superposition of diffusion processes over a continuous range $(0,c^*]$ of 
speeds does not satisfy the strongly continuous semigroup property, but the splitting up 
of masses take place at each time point $t>\tau^*$ for some $\tau^*>0$.} 
As a consequence, not each space-time point has a unique velocity vector, but 
the \emph{general continuity equation} 
\begin{equation}\label{gconteq}
\begin{aligned}
    \frac{\partial v}{\partial t}(\x,t) + \nabla \cdot \vj(\x,t) = 0  
   \qquad\quad  \mbox{for} \qquad \quad  t\not=m\,\tau, \quad m\in\N_0
\end{aligned}
\end{equation}
holds with \emph{flux density} 
\begin{equation}\label{defj}
\begin{aligned}
      \vj(\x,t) =
       - c_0\,\int_{S_1(\vO)}  \frac{v(\x+R(t)\,\y,n(t)\,\tau)}{|S_1(\vO)|}\,\y\,\d \sigma(\y)\,.
\end{aligned}
\end{equation} 
For those coordinates $(\x,t)$ for which
$$
   v(\x + R(t)\,\y,n(t)\,\tau) 
         \approx  v(\x ,n(t)\,\tau) 
                + R(t)\,\nabla v(\x ,n(t)\,\tau)\cdot \y\,
$$
holds, it follows that 
\begin{equation}\label{Dt}
\begin{aligned}
      \vj(\x,t) \approx 
        - D(t)\, \nabla v(\x,n(t)\,\tau)  \qquad\mbox{with}\qquad D(t):=\frac{c_0\,R(t)}{N}\,,
\end{aligned}
\end{equation} 
since
\begin{equation}\label{Dkl}
\begin{aligned}
  \int_{S_1(\vO)}  y_k\,y_l \,\d \sigma(\y)
       = \delta_{k,l}\,\, |S_1(\vO)|/N \qquad\quad k,\,l\in\{1,2,\ldots,N\}\,.
\end{aligned}
\end{equation} 
But this approximation of the flux density resembles \emph{Fick's law}\footnote{Cf. 
Remark~\ref{rema:limit} at the end of Section~\ref{sec-causalmodel}.}
$$
   \vj(\x,t) = - D_0\nabla v(\x,t)  \qquad\mbox{with}\qquad  
   D_0:=D(\tau/2)\,.
$$
Note if $\tau\to 0$ say under the side condition $D_0(\tau/2)=\const.$, 
then $c_0\to \infty$, i.e. causality does not hold for this limit.

Finally, we show that the diffusion process~(\ref{defcausdiff10}) satisfies an inhomogeneous wave 
equation with time dependent coefficient such that mass is conserved.  \\

This paper is organized as follows: In Section~\ref{sec-semigroups} we prove under fairly 
reasonable assumptions that causal diffusion cannot satisfy the property of a strongly 
continuous semigroup. 
For the convenience of the reader we put some of the technical parts of the proof of the 
main theorem in the Appendix. 
Subsequently we present and analyze a causal diffusion model in Section~\ref{sec-causalmodel}. 
The paper concludes with a short section of conclusion (Section~\ref{sec-conclusion}).

\section{Real-valued semigroups and causality}
\label{sec-semigroups}

It is common in science to assume that the operator $\A$ in~(\ref{evolequv}) has nice 
properties like linearity, translation invariance and some kind of smoothness such that 
it can infered that $\A$ is a space convolution operator. Theorem 4.2.1 
in~\cite{Hoe03} provides such a list of necessary properties. For example, translation and 
rotation invariance of $\A$ means that the medium in which the process takes place is 
homogeneous and isotropic. In this section we assume 
that $\A$ is a \emph{space convolution operator} and an \emph{(infinitesimal) generator} 
of a strongly continuous semigroup in $L^1(\R^N)$ and prove that the solution of evolution 
equation~(\ref{evolequv}) cannot satisfy causality condition~(\ref{causaGcond}). 
This result suggests that even if some of the nice properties of $\A$ are dropped, causality cannot 
be restored. And therefore we present a new model for causal diffusion that do not satisfy 
the semigroup property for all time points in the following section.

For the convenience of the reader we put the technical parts of the proof of the main theorem 
in the Appendix. \\

\subsection{Description of the evolution problem}

The notion of the Fourier transform and the convolution used in this paper are specified 
at the beginning of the Appendix. Moreover, we use the following notion

\begin{defi}\label{defaA}
Let $N\in\N$ and $\hat\a:\R^N\to\R^N$ be a positive and rotational symmetric function such that
\begin{equation*}
\begin{aligned}
   \exp(-\hat \a(\k)) \in \S(\R^N)
\qquad\mbox{and}\qquad
   (2\,\pi)^{-N/2}\,\int_{\R^N} |\F^{-1}\{\exp(-\hat \a(\k))\}(\x)|\,\d \x=1\,.
\end{aligned}
\end{equation*}
Moreover, let
$\a := (2\,\pi)^{-N/2}\,\F^{-1}\{\hat \a\}$.
The operator $\A:\S'(\R^N)\to \S'(\R^N)$ is defined as the convolution operator
\begin{equation*}
\begin{aligned}
   (\A\,u)(\x) := -(\a *_\x u)(\x) \qquad\quad \mbox{ for }\qquad u\in \S'(\R^N)\,.
\end{aligned}
\end{equation*}
Here $*_\x$ denotes the convolution with respect to the space variable $\x\in\R^N$.
\end{defi}

This definition is reasonable, since the convolution $f*_\x g$ of distributions $f\in\S(\R^N)$
and $g\in\S'(\R^N)$ is well-defined and lies in $\S'(\R^N)$ (cf. e.g.~\cite{GasWit99,Hoe03}).\\

Let $\hat\a$ and $\A$ be as in Definition~\ref{defaA}. Then
\begin{equation}\label{defG}
\begin{aligned}
   G(\x,t) :=  (2\,\pi)^{-N/2}\,\F^{-1}\{\exp(-\hat\a\,t)\}(x) 
\qquad\quad \mbox{($\x\in\R,\,t\in\R$)}\,
\end{aligned}
\end{equation}
is well-defined and defines by
\begin{equation}\label{relStG}
\begin{aligned}
   (S_t\,u)(\x) = (G(\cdot,t) *_\x u)(x) \qquad\quad \mbox{ for $u\in \S'(\R)$}\,
\end{aligned}
\end{equation}
a semigroup $\{S_t:\S'(\R^N)\to \S'(\R^N)\,|\,t\geq 0\}$ with generator $\A$ (cf.
e.g. Chapter XVII in~\cite{DauLio92_5} or Chapter 7.4 in~\cite{EVa99}).
Moreover, $G$ satisfies the evolution equation
\begin{equation}\label{evoleqG}
\begin{aligned}
     \frac{\partial G}{\partial t}  - \A\,G = 0
               \qquad \mbox{ on }\qquad \R^N\times (0,\infty)
\end{aligned}
\end{equation}
with initial condition
\begin{equation}\label{initevoleqG}
\begin{aligned}
    G(\x,0) = \delta(\x)\,.
\end{aligned}
\end{equation}
Here $\delta(\x)$ denotes the \emph{dirac delta distribution} on $\R^N$. 
Conversely, if $\{ S_t \,|\, t \geq 0 \}$ is a semigroup with a generator $\A$ that is a 
space convolution operator with kernel $a$, then $G$ defined as in~(\ref{defG}) 
satisfies~(\ref{relStG}) and~(\ref{evoleqG}) with~(\ref{initevoleqG}).

\begin{rema}
The assumption that $\hat\a$ is rotational symmetric is motivated by the following fact. If a 
process takes place in a homogeneous and isotropic medium, then - by symmetry - the fundamental 
solution $G_0$ of the respectively evolution equation must be rotational symmetric in $\x$. 
And consequently, $\hat\a$ is rotational symmetric, too. Since this paper is a first approach on
causality of semigroups and its evolution equations, we focus on these special types of processes. 
\end{rema}

\begin{exam}\label{ex:gauss}
A prominent semigroup is defined by the \emph{Gauss function}
\begin{equation}\label{defGauss}
\begin{aligned}
   G(\x,t) = (4\,\pi\,D_0\,t)^{-n/2}\,\exp\left(-\frac{|\x|^2}{4\,D\,t} \right)
         \qquad (\x\in\R^N,\,t\in\R)\,,
\end{aligned}
\end{equation}
where $D_0$ is a  constant.
The respective generator $\A$ is defined by
\begin{equation*}
\begin{aligned}
   \hat \a(\k) := D_0\,|\k|^2\, \qquad\quad (\k\in\R^N)\,.
\end{aligned}
\end{equation*}
The process described by this semigroup is not ``causal'', since
the Gauss function $G$ is everywhere positive. Moreover, $G$ satisfies the standard diffusion 
equation~(\ref{stddiffeq}) with $u(\x)=\delta(\x)$. 
$\Box$
\end{exam}

\subsection{The structure of $\hat a$}

We first investigate the properties of $\hat a$, the Fourier transform of the kernel $a$ 
of the operator $\A$. For this purpose we recall some facts about \emph{holomorphic roots.}

If $f$ and $g$ are entire functions such that $f=g^n$, then we call
$g$ an $n-$th holomorphic root of $f$. If the entire function $f$
has an $n-$th holomorphic root and
\begin{equation}\label{roots}
\begin{aligned}
  f^{(0)}(c)=0,\quad  f^{(1)}(c)=0,\quad\ldots,\quad  f^{(m-1)}(c)=0,\quad
        f^{(m)}(c)\not=0\,.
\end{aligned}
\end{equation}
for some $c\in\C$, then $n$ divides $m$. Here $f^{(0)}:=f$ and $f^{(n)}$ denotes
the $n-$th derivative of $f(z)$ with respect to $z$. For a proof of this statement see the
``Wurzelkriterium'' in Chapter 3 Section 1 in~\cite{Rem95b}.

\begin{theo}\label{holroot}
Let $f$ be an entire function. If $f$ has no zeros, then $f$ has an
$n-$th holomorphic root for each $n\in\N$. Conversely, if $f$ has an
$n-$th holomorphic root for each $n\in\N$, then $f$ has no zeros.
\end{theo}

\begin{proof}
The first statement is well-known (cf. for example the ``Wurzelsatz'' in
Chapter 3 Section 2 in~\cite{Rem92a}).
For the second statement. Assume that $f$ has a holomorphic root for each
$n\in\N$. Then for each $c\in\C$ there exists a finite
$m\in\N\cup\{0\}$ such that~(\ref{roots}) holds. According to the previous
remark, each $n\in\N$ divides $m$ and consequently $m$ must be zero.
But this means that $f$ has no zeros which concludes the proof.
\end{proof}

Now, with the help of Theorem~\ref{holroot} and the Paley-Wiener-Schwartz Theorem, 
we can prove that $z\in\C^N\mapsto \hat a(z)$ is entire if~(\ref{causaGcond}) holds.

\begin{lemm}\label{lemm:1}
If $G$ defined by~(\ref{defG}) satisfies condition~(\ref{causaGcond}),
then $\hat\a$ can be extended to an entire function denoted by
$z\in\C\mapsto \hat\a(z)$.
\end{lemm}

\begin{proof} 
Let
\begin{equation}\label{defft}
\begin{aligned}
  f_t(\k):= \F\{G(\cdot,t)\}(\k)  \qquad\mbox{for}\qquad k\in\R^N,\,t>0\,.
\end{aligned}
\end{equation}
If $G$ satisfies~(\ref{causaGcond}), then according to the Paley-Wiener-Schwartz
Theorem (cf. Theorem~\ref{th:PWS1} in the Appendix),
each $f_t(\k)$ can be extended to an entire function, denoted by $F_t(\z)$, such that
\begin{equation}\label{Ftft}
\begin{aligned}
  F_t(\k) = f_t(\k)=\exp\{-\hat\a(\k)\,t\}  \qquad\mbox{for}\qquad \k\in\R^N\,.
\end{aligned}
\end{equation}
Moreover, it follows that $F_{1/n}(\z)$ is an $n-$th holomorphic root of $F_1(\z)$, i.e.
$(F_{1/n}(\z))^n=F_1(\z)$ and $F_{1/n}(\z)$ is entire.
This together with Theorem~\ref{holroot}, implies that $F_1$ has no zeros.
Hence there exists an entire function $\hat\a_0(\z)$
such that
\begin{equation*}
\begin{aligned}
   F_t(\z) = F_t(\vO)\,\exp\{-\hat\a_0(\z)\,t\}  \qquad\mbox{with}\qquad
    \hat\a_0(\z) = -\int_{\gamma_z} \frac{F_t'(\xi)}{F_t(\xi)}\,\d \xi\,,
\end{aligned}
\end{equation*}
where $\gamma_z$ is any path connecting $0$ and $\z$. Since~(\ref{Ftft}) holds,
it follows that $\hat\a_0(\k)=\hat\a(\k)$ for $\k\in\R$ and
$F_t(\vO)=1$. This proves the lemma.
\end{proof}

\subsection{Proof of noncausality}

For the convenience of the reader we recall some notion about polynomials in several variables.
Let $z=(z_1,\ldots,z_N)\in\C^N$ and $n:=(n_1,\ldots,\,n_N)\in\N^N$. Then
\begin{equation}\label{polynomial}
    z^n := z_1^{n_1}\,\ldots\,z_N^{n_N}  \qquad \mbox{and}\qquad
  |z|^2:=\sum_{n=1}^N z_n^2\,.
\end{equation}
In the following $ch\,\supp(g)$ denotes the \emph{closed convex hull} of the support of $g$.

\begin{theo}\label{theo:main}
Let $\hat\a$ and $\A$ be defined as in Definition~\ref{defaA} and $G$ be defined as
in~(\ref{defG}). Moreover, let $f_1:=G(\cdot,1)$. If $\hat\a$ is entire and 
\begin{equation}\label{A2} 
   ch\,\supp(f_1) = B_{c_0}(\vO)  \qquad \mbox{with}\qquad c_0=1\,, 
\end{equation}
then $\supp(\A f_1)\subseteq B_{c_0}(\vO)$ cannot hold.
\end{theo}

\begin{proof}
We perform a proof by contradiction and assume $\supp(\A f_1)\subseteq  B_1(\vO)$.
There are two cases, either $\F^{-1}(\hat a)$ has or has not compact support.
\begin{itemize}
\item [i)] Assume that $a$ has compact support.
      Then according to the support theorem for distributions with compact support
      (cf. Theorem~4.3.3 in~\cite{Hoe03}), we have
      \begin{equation*}
      \begin{aligned}
       B_1(\vO)
         \supseteq ch\,\supp (\A f_1)
          &= ch\,\supp(\F^{-1}(\hat a)) + ch\,\supp(f_1) \\
          &= ch\,\supp(\F^{-1}(\hat a)) + B_1(\vO)\,,
      \end{aligned}
      \end{equation*}
      since $\A f_1 = \F^{-1}(\hat a) *_\x f_1$.
      Only if $\F^{-1}(\hat a)$ has singular support $\{\vO\}$, we do not obtain a contradiction.
      If $\supp(\F^{-1}(\hat a))=\{\vO\}$, then $\hat \a(z)$ is a polynomial in
      $\i\,z=\i\,z_1\ldots\,z_N$ (cf. Theorem 2.3.4 in~\cite{Hoe03}).
\begin{itemize}
\item [a)] If $\hat\a(z)=b\in\R$, then $f_1(z)=\exp\{-b\}$ and thus
$\supp(f_1)=\{\vO\}\not=B_1(\vO)$ (cf.~(\ref{A2})).

\item [b)] If $\hat\a(z)=\i\,b\,z$ with $b\in\R\backslash\{0\}$, then $f_1(z)=\exp\{-\i\,b\,z\}$
and thus $\supp(f_1)=\{(b,\ldots,b)\}\not=B_1(\vO)$. (For $b=0$ we have case i) a).)
\item [c)]
In Theorem~\ref{caseP} in the Appendix, we show that if $\hat \a(z)$ is a polynomial in $\i\,z$ 
with exponent larger than $1$, then $\supp(f_1)$ cannot have compact support.
\end{itemize}
For each of the cases a)-c) we obtain a contradiction which proves that the assumption
$\supp(\A f_1)\subseteq B_1(\vO)$ cannot hold for case i).

\item [ii)] That $\supp(\A f_1)\subseteq B_1(\vO)$ cannot hold if $a$ does not have compact 
            support, follows from an 
            estimation of the minimum modulus of $\hat a$. The details are carried out 
            in Theorem~\ref{theo:caseSeries} in the Appendix.
\end{itemize}
In summary, we have shown that $\supp(\A f_1)\subseteq B_1(\vO)$ cannot hold for both cases. 
\end{proof}

Now we are ready to prove the noncausality of $G$. More precisely 

\begin{theo}\label{theo:main0}
Let $N\in\N$ and $\hat\a\in\S(\R^N)$ be a function that is real valued and rotational 
symmetric. Then $G$ defined by~(\ref{defG}) does not satisfy causality 
condition~(\ref{causaGcond}).
\end{theo}

\begin{proof}
We perfom a proof by contradiction. Assume that
\begin{itemize}
\item [(AS)] $\quad G$ satisfies condition~(\ref{causaGcond})
\end{itemize}
and without loss of generality, we can assume that~(\ref{A2}) holds,   
since $\F^{-1}\{\hat\a\}(\x)$ is rotational symmetric. (This can be accomplished by rescaling of 
space and time.) 
Since $G$ defined by~(\ref{defG}) solves evolution equation~(\ref{evoleqG}),
it follows
$$
     ch\,\supp (\A\,G(\cdot,t))
           =  ch\,\supp \left(\frac{\partial G}{\partial t}(\cdot,t)\right)
            \subseteq  ch\,\supp (G(\cdot,t))
            \subseteq B_{t}(\vO)
$$
for each $t>0$. In particular,
\begin{equation}\label{causpropA}
\begin{aligned}
          ch\,\supp(\A\,G(\cdot,1))\subseteq B_1(\vO)\,.
\end{aligned}
\end{equation} 
But according to Theorem~\ref{theo:main}, property~(\ref{causpropA}) cannot hold.
This contradiction shows that the assumption (AS) is not true and concludes the proof.
\end{proof}

\section{A causal model for diffusion}
\label{sec-causalmodel}

As shown in Theorem~\ref{theo:main0}, it is not possible that causal diffusion satisfies the
semigroup property for every time point $t>0$. In this section we present a model for
diffusion that satisfies
\begin{itemize}
\item the semigroup property for a discrete set of time points $\{\tau_n\,|\,n\in\N_0\}$,

\item an inhomogeneous wave equation with a time dependent coefficient such that

\item the total mass is conserved and

\item causality holds.
\end{itemize}

\vspace{0.4cm}
\noindent
In the following we use the notion
\begin{itemize}
\item [A1)] $c,\,\tau\in (0,\infty)$ and $\tau_m:=m\,\tau$ for $m\in\N_0$.

\item [A2)] If $t\in (\tau_m,\tau_{m+1}]$ ($m\in\N_0$), then $n(t):=m$ and 
            $R(t):=c\,(t-\tau_{n(t)})$.

\item [A3)] $f:[0,c^*]\to [0,\infty)$ is such that $\int_0^{c^*} f(c) \,\d c=1$.

\end{itemize}
We emphasize that if $t=\tau_k$, then $n(t)=k-1$ and $R(t) = c_0\,\tau$. 
Moreover, $\d \sigma(\x')$ denotes the Lebesgue surface measure on $\R^N$ and
$|S_{R}(\vO)|$ denotes the surface area of the sphere $S_R(\vO)\subset \R^N$.
We note that
\begin{equation}\label{propS}
\begin{aligned}
      &\int_{S_R(\x)} f(\x')\,\d \sigma(\x') = \int_{S_1(\vO)} f(\x+R\,\y)\,R^{N-1}\, \d \sigma(\y)
\qquad\mbox{and} \\
    &\qquad\qquad\qquad |S_{R}(\vO)|=|S_1(\vO)|\,R^{N-1} \,.
\end{aligned}
\end{equation}

\subsection{Definition of the model}

First we consider the case where each particular process takes place with the
same speed $c$, i.e. $c$ and $\tau$ are fixed in time.
\begin{defi}\label{def:diff1}
Let $c$, $\tau$, $\tau_m$, $n(t)$ and $R(t)$ be as in A1)-A2). The concentration $v$ of a
substance  diffusing via a constant speed $c$ with initial concentration $u\in L^1(\R^N)$
is defined by
\begin{equation}\label{defcausdiff1}
\begin{aligned}
   v(\x,t)
      =\int_{S_{R(t)}(\x)} \frac{ v\left(\x',\tau_{n(t)}\right)
       }{|S_{R(t)}(\vO)|} \,\d \sigma(\x')\,
\quad\mbox{with}\quad v|_{t=0}=u\,.
\end{aligned}
\end{equation} 
In analogy to wave mechanics, we call $\lambda:=c\,\tau$ and 
$k:=(2\,\pi)/\lambda$ the \emph{wave length} and the \emph{wave number} of the diffusion 
process.
\end{defi}

This definition means that the concentration on the sphere $S_{R(t)}(\x)$ 
at time $\tau_{n(t)}$ determines the concentration in point $\x$ at time $t$.  
If we consider the paths of the particles as lines during the time period 
$(\tau_{n(t)},\tau_{n(t)+1})$, then the particle speed is 
$$
    \frac{R(t)}{t-\tau_{n(t)}} = c\,.
$$
A more detailed interpretation of this model is given in Subsection~\ref{subsec-interpr}.

\begin{rema}
For the one dimensional case $N=1$, we have $S_{R(t)}(x) = \{x-R(t),x+R(t)\}$,
$|S_{R(t)}(x)|=2$ and
\begin{equation*}
\begin{aligned}
   \int_{S_{R(t)}(x)}  \d \sigma(x')
      \equiv \int_\R ( \delta(x'-R(t)) + \delta(x'+R(t)) )\,\d x'\,.
\end{aligned}
\end{equation*}
Therefore Definition~(\ref{defcausdiff1}) reads as follows
\begin{equation*}
\begin{aligned}
   v(x,t)
      = \frac{1}{2}\,( v(x-R(t),\tau_{n(t)}) + v(x+R(t),\tau_{n(t)}) )\,.
\end{aligned}
\end{equation*}
For this special case $x\mapsto v(x,t)$ for each $t>0$ is a discrete measure on $\R$.
\end{rema}

That mass is conserved by the diffusion model~(\ref{defcausdiff1}) follows from
the following lemma.

\begin{lemm}\label{lemm:masscons}
Let $u\in L^1(\R^N)$ and $v$ be defined as in~(\ref{defcausdiff1}).
Then $v$ is well-defined and satisfies
$$
    \int_{\R^N} v(\x,t)\,\d \x =\int_{\R^N} u(\x)\,\d \x
  \qquad \mbox{for}  \qquad   t>0\,.
$$
\end{lemm}

\begin{proof}
We perform a proof by induction.  \\
Let $n=0$, i.e. $t\in (0,\tau]=(\tau_0,\tau_1]$. Then~(\ref{defcausdiff1}) and~(\ref{propS})
imply
\begin{equation}\label{repv01}
\begin{aligned}
   v(\x,t)
      = \int_{S_1(\vO)}  \frac{u(\x+R(t)\,\y)}{|S_1(\vO)|}
         \,\d \sigma(\y)
     \qquad\quad (t\in (0,\tau])\,
\end{aligned}
\end{equation}
which shows that $v(\cdot,t)$ is well defined for $t\in (0,\tau]$ and satisfies
\begin{equation*}
\begin{aligned}
   \int_{\R^N} v(\x,t)\,\d \x
      &= \int_{S_1(\vO)} \int_{\R^N} u(\x+R(t)\,\y) \,\d \x\, \frac{1}{|S_1(\vO)|}
        \,\d \sigma(\y) \\
      &= 1\,  \int_{\R^N} u(\x) \,\d \x\,.
\end{aligned}
\end{equation*}
Now let $n\in\N$, i.e. $t\in (\tau_n,\tau_{n+1}]$. Analogously,~(\ref{defcausdiff1})
and~(\ref{propS}) imply
\begin{equation}\label{repv02}
\begin{aligned}
   v(\x,t)
      = \int_{S_1(\vO)}  \frac{v(\x+R(t)\,\y,\tau_n)}{|S_1(\vO)|}
         \,\d \sigma(\y)\,
     \qquad\quad (t\in (\tau_n,\tau_{n+1}] )\,
\end{aligned}
\end{equation}
which shows by induction that $v(\cdot,t)$ is well defined for $t\in (\tau_n,\tau_{n+1}]$.
Moreover, from the last identity, it follows that
\begin{equation*}
\begin{aligned}
   \int_{\R^N} v(\x,t)\,\d \x
      &= \int_{S_1(\vO)} \int_{\R^N} v(\x+R(t)\,\y,\tau_n)\,\d \x\, \frac{1}{|S_1(\vO)|}
         \,\d \sigma(\y) \\
      &= 1\,  \int_{\R^N} v(\x,\tau_n) \,\d \x\,
\end{aligned}
\end{equation*}
which together with the induction assumption
$\int_{\R^N} v(\x,t) \,\d \x=\int_{\R^N} u(\x) \,\d \x$ for $t\in (0,\tau_n]$ imply that
$$
   \int_{\R^N} v(\x,t)\,\d \x = \int_{\R^N} v(\x,\tau_n)\,\d \x =\int_{\R^N} u(\x) \,\d \x\,.
$$
\end{proof}

Now we define diffusion consisting of various noninterfering processes that take place
with different speeds $c\in (0,c^*]$.  

\begin{defi}\label{def:diff2}
Let $c^*,\tau^*>0$, $c\in (0,c^*]$ and $\tau\in [\tau^*,\infty)$ be such that $c\,\tau=c^*\,\tau^*$, 
$f$ be as in A3)  and $u\in L^1(\R^N)$ be such that a solution $u_{c}$ 
of\footnote{This problem is ill-posed and thus cannot be solved without additional information.}
\begin{equation}\label{defcausdiff2b}
\begin{aligned}
   u(\x) = \int_0^{c^*} f(c)\,u_{c}(\x) \,\d c\,
         \qquad\quad \mbox{($c\,\tau=const$)}\,
\end{aligned}
\end{equation}
exists.
Moreover, let $v_{c}$ be defined as in~(\ref{defcausdiff1}) with initial data $u_{c}$.
Diffusion of an initial concentration $u\in L^1(\R^N)$ with speed distribution $f$
is defined by
\begin{equation}\label{defcausdiff2}
\begin{aligned}
   v(\x,t) = \int_0^{c^*} f(c)\,v_{c}(\x,t) \,\d c
\qquad\mbox{with}\qquad v|_{t=0}=u\,.
\end{aligned}
\end{equation}
We call $\lambda:=c^*\,\tau^*$ and $k:=(2\,\pi)/\lambda$ the \emph{wave length} and 
the \emph{wave number} of the diffusion process. 
\end{defi}

That model~(\ref{defcausdiff2}) conserves mass follows at once from
Lemma~\ref{lemm:masscons} and property $\int_0^{c^*} f(c) \,\d c=1$.

\subsection{The Green function of diffusion and its properties}

\begin{defi}\label{defi:Green}
Let $p(u;\x,t)$ ($\x\in\R^N$, $t\geq 0$) denote the concentration of the diffusion processes 
with initial distribution $u\in L^1(\R^N)$, i.e. $p(u;\x,0)=u$. 
If
$$
           p(u;\x,t) = p(\delta(\x);\x,t) *_\x u
   \qquad\mbox{for each}\qquad u\in L^1(\R^N)\,,
$$
then we call $G(\x,t):=p(\delta(\x);\x,t)$ the \emph{Green function} of diffusion.
\end{defi}

Below in Proposition~\ref{prop:green02} and Remark~\ref{rema:green01} we show that the
Green function exists if and only if $f(c) =\delta(c-c_0)$.
But before we can do that we need a proposition.

\begin{theo}\label{theo:propG}
Let $\tau$, $c$, $n(t)$ and $R(t)$ be as in A1)-A2) and $G$ be defined by~(\ref{defcausdiff1})
with $u=\delta(\x)$. Then $G$ satisfies\footnote{This is a 
modification of the semigroup property of continuous model processes in Stochastic Analysis.}
\begin{equation}\label{mainpropG}
\begin{aligned}
   G(\cdot,t)
     = G(\cdot,\tau_{n(t)}) *_\x G(\cdot,t-\tau_{n(t)})
        \qquad \mbox{for}\qquad t> 0\,
\end{aligned}
\end{equation}
and 
\begin{equation*}
\begin{aligned}
   G(\cdot,\tau_{k+m})
     = G(\cdot,\tau_k) *_\x G(\cdot,\tau_m)
       \qquad \mbox{for}\qquad k,\,m\in\N_0\,.
\end{aligned}
\end{equation*}
\end{theo}

\begin{proof}
By~(\ref{propS}) we have for $R_0>0$, $\x\in\R^N$ and $g\in L^1(\R^N)$:
\begin{equation*}
\begin{aligned}
   \int_{\R^N} f(\x')\,&\delta(R_0-|\x-\x'|) \,\d \x'\\
         &= \int_0^\infty \int_{S_1(\vO)}  f(\x+R\,\y)\,\delta(R_0-|R\,\y|) \,R^{N-1}
                                       \,\d \sigma(\y)\,\d R\\
         &= \int_{S_1(\vO)}  f(\x+R_0\,\y)\,R_0^{N-1} \,\d \sigma(\y)
         = \int_{S_{R_0}(\x)}  f(\x')\, \,\d \sigma(\x')\,.
\end{aligned}
\end{equation*}
From this and~(\ref{defcausdiff1}), it follows that
\begin{equation*}
\begin{aligned}
   G(\x,t)
     &= \frac{
              \int_{\R^N} G(\x',\tau_n) \, \delta(R(t) - |\x-\x'|) \,\d \x'
             }{|S_{R(t)}(\vO)|}\,.
\end{aligned}
\end{equation*}
In particular, since $G(\x',0)=\delta(\x')$,
\begin{equation}\label{repG1}
\begin{aligned}
   G(\x,t)
     &= \frac{
              \delta(c_0\,t - |\x|)
             }{|S_{c_0\,t}(\vO)|}\,  \qquad \mbox{for}\qquad t\in(0,\tau]\,.
\end{aligned}
\end{equation}
Combining the last two results yields
\begin{equation*}
\begin{aligned}
   G(\x,t)
     &=  \int_{\R^N} G(\x',\tau_n) \, G(\x-\x',t-\tau_n) \,\d \x'\, \\
     &= (G(\cdot,\tau_n) *_\x G(\cdot,t-\tau_n))(\x)
        \qquad \mbox{for}\qquad t>0,\,\,x\in\R^N\,.
\end{aligned}
\end{equation*}
Moreover, we get from this for $t=\tau_{n+1}$ that
\begin{equation}\label{repG2}
\begin{aligned}
   G(\x,\tau_{n+1})
     = (G(\cdot,\tau_n) *_\x G(\cdot,\tau))(\x) \,
\end{aligned}
\end{equation}
and consequently by induction
$G(\cdot,\tau_{k+m}) = G(\cdot,\tau_k) *_\x G(\cdot,\tau_m)$. As was to be shown.\\
\end{proof}

\begin{prop}\label{prop:green02}
Let $v$ be as in Definition~\ref{def:diff1} and $G$ be defined by~(\ref{defcausdiff1})
with $u(\x)=\delta(\x)$. Then
\begin{equation*}
\begin{aligned}
   v(\x,t) = (G(\cdot,t) *_\x u)(\x)
   \quad\qquad t>0,\,\x\in\R^N\,,
\end{aligned}
\end{equation*}
i.e. $G$ is the Green function of diffusion with constant speed $c$. 
In particular, $v$ satisfies\footnote{This equation is a modification of the 
Chapman-Kolmogorov equation and since $G(\cdot,t)$ is a positive distribution, 
$\mu_t(\x):=G(\x,t)\,\d\x$ is a measure on $\R^N$.}
\begin{equation*}
\begin{aligned}
   v(\x,t)
     = \int_{\R^N} G(\x',\tau_{n(t)}) \, v(\x-\x',t-\tau_{n(t)})\,\d \x'
   \quad\qquad t>0,\,\x\in\R^N\,
\end{aligned}
\end{equation*}
with initial condition $v(\x,0)=u(\x)$.
\end{prop}

\begin{proof}
We prove the first statement of the proposition by induction. \\
Let $n=0$, i.e. $t\in (0,\tau]=(\tau_0,\tau_1]$. According to~(\ref{defcausdiff1})
and~(\ref{propS}), we have
\begin{equation}\label{repG1b} 
\begin{aligned}
   G(\x,t)
      =  \int_{S_1(\vO)} \frac{\delta(\x+R(t)\,\y)}{|S_1(\vO)|} \, \d \sigma(\y)
\qquad (t\in (0,\tau])\,.
\end{aligned}
\end{equation}
From this and~(\ref{repv01}), we infer for $t\in (0,\tau]$:
\begin{equation*}
\begin{aligned}
   v(\x,t)
      &= \int_{S_1(\vO)}\frac{ u(\x+R(t)\,\y)
         }{|S_1(\vO)|}  \, \d \sigma(\y)\, \\
      &= \int_{\R^N} u(\x') \, \int_{S_1(\vO)}   \frac{ \delta(\x+R(t)\,\y-\x')
         }{|S_1(\vO)|} \, \d \sigma(\y) \,\d \x' \\
     &= ( u *_\x G(\cdot,t) )(\x) \,.
\end{aligned}
\end{equation*}
This proves case $n=0$. \\
Assume that $v(\x,t) = ( G(\cdot,t) *_\x u)(\x)$ holds for
$t\in (\tau_{n-1},\tau_n]$ (induction assumption). Now let $t\in (\tau_n,\tau_{n+1}]$.
Analogously as above, it follows from~(\ref{repv01}) that
\begin{equation*}
\begin{aligned}
   v(\x,t)
      &= \int_{S_1(\vO)}\frac{ v(\x+R(t)\,\y,\tau_{n(t)})
         }{|S_1(\vO)|}  \, \d \sigma(\y)\, \\
      &= \int_{\R^N} v(\x',\tau_{n(t)}) \, \int_{S_1(\vO)}   \frac{ \delta(\x+R(t)\,\y-\x')
         }{|S_1(\vO)|} \, \d \sigma(\y) \,\d \x' \\
     &= ( v(\cdot,\tau_{n(t)}) *_\x G(\cdot,t-\tau_{n(t)}) )(\x) \,
\end{aligned}
\end{equation*}
which together with the induction assumption and~(\ref{mainpropG}) imply
\begin{equation*}
\begin{aligned}
   v(\cdot,t)
     = [ G(\cdot,\tau_n) *_\x u ] *_\x G(\cdot,t-\tau_n) 
     = G(\cdot,t) *_\x u\,.
\end{aligned}
\end{equation*}
This proves the first claim of the proposition. \\
According to~(\ref{mainpropG}) and the first statement of the proposition, we have
$$
    v(\cdot,t) 
           = G(\cdot,t) *_\x u 
           = G(\cdot,t-\tau_n) *_\x G(\cdot,\tau_n) *_\x u
           = G(\cdot,\tau_n) *_\x v(\cdot,t-\tau_n)
$$
which proves the last claim of the proposition.
\end{proof}

\begin{rema}\label{rema:green01}
Let $v_c$ be as in Definition~\ref{def:diff1} and $v$, $f$ be as in Definition~\ref{def:diff2}.
Moreover, let $G_c$ be defined by~(\ref{defcausdiff1}) with $u_c=\delta(\x)$ and
$$
     \tilde G(\cdot,t) := \int_0^{c^*}  f(c)\,G_c\left(\cdot,t\,\right) \,\d c \,.
$$
From $\int_0^{c^*} f(c) \,\d c=1$ and $G_c(\x,0) = \delta(\x)$, we infer 
$\tilde G(\x,0) = \delta(\x)$. However, in general (cf. Definition~\ref{def:diff2} and 
Definition~\ref{defi:Green})
\begin{equation*}
\begin{aligned}
   v(\x,t) &= \int_0^{c^*}  f(c)\, (G_c\left(\cdot,t\,\right)*_\x u_c)(\x) \,\d c \\
           &\not= \int_0^{c^*}  f(c)\, (G_c\left(\cdot,t\,\right)*_\x u)(\x) \,\d c
   = (\tilde G\left(\cdot,t\,\right)*_\x u)(\x)
\end{aligned}
\end{equation*}
and therefore we do not call $\tilde G$ Green function if $f(c)\not=\delta(c-c_0)$.
\end{rema}

The following two propositions are about the smoothness properties of the diffusion process 
defined as in Definition~\ref{def:diff1} and Definition~\ref{def:diff2}, respectively. 
In the following, we use the subscript $c$ only if necessary.

\begin{prop}\label{prop:vcont}
Let $v$ be defined as in Definition~\ref{def:diff1} with $u\in L^1(\R^N)$. \\
a) If $u$ is continuous, then $v(\x,\cdot): \R\to (0,\infty)$ is continuous for 
each $\x\in\R^N$ and 
$$
          v(\x,0) := v(\x,0+) = u(\x)\,.
$$
b) If $u$ is differentiable\footnote{We can also assume Gateaux differentiabilty.}, 
then $\frac{\partial v}{\partial t} (\x,\cdot): \R\to (0,\infty)$ 
is continuous on $(\tau_m,\tau_{m+1})$ for each $m\in\N$ and 
\begin{equation}\label{jumpcond}
\begin{aligned}
   \left[\frac{\partial v}{\partial t}\right]_{t=\tau_m-}^{t=\tau_m+}
          =  -u *_\x \frac{\partial G}{\partial t}(\cdot,\tau_m-) 
           \qquad\mbox{for} \qquad m\in\N_0\,.
\end{aligned}
\end{equation}
In particular,
$$
   \frac{\partial v}{\partial t}(\cdot,\tau_m) 
          :=  \frac{\partial v}{\partial t}(\cdot,\tau_m+) = 0 
           \qquad\mbox{for} \qquad m\in\N_0\,.
$$
\end{prop}

\begin{proof}
a) Let $u$ be continuous. 
From~(\ref{repv01}), it follows that 
\begin{equation*}
\begin{aligned}
   &| v(\x,t+s)- v(\x,t)|   \\
     &\qquad \leq   \int_{S_1(\vO)}  \frac{| u(\x+R(t+s)\,\y)- u(\x+R(t)\,\y) |}{|S_1(\vO)|}  
                     \,\d \sigma(\y)
\end{aligned}
\end{equation*}
and because $R(t)$ is continuous for $t\in (0,\tau)$, it follows that $t\mapsto v(\x,t)$  is 
continuous for $t\in (0,\tau)$.  
Let $t=\tau$ and $s\in (-\tau,0)$. Then we have (cf. A2))
$$
    R(t) = c_0\,\tau   \qquad\mbox{and}\qquad   R(t+s) = c_0\,(\tau-|s|)
$$
and the previous estimation reads as follows 
\begin{equation*}
\begin{aligned}
   &| v(\x,t+s)- v(\x,t)|    \\
     &\qquad \leq   \int_{S_1(\vO)} \frac{| u(\x+c_0\,(\tau-|s|)\,\y)- u(\x+c_0\,\tau\,\y) | }
                {|S_1(\vO)|} \,\d \sigma(\y)
\end{aligned}
\end{equation*}
which shows that $t\mapsto v(\x,t)$ is continuous from the left at $t=\tau$.  \\
Let $t=0$ and $s\in (0,\tau)$, then~(\ref{defcausdiff1}) implies 
\begin{equation*}
\begin{aligned}
   v(\x,0+)  
     &= \lim_{s\to 0,\,s>0}\int_{S_1(\vO)}   \frac{  u(\x +R(s)\,\y) }{|S_1(\vO)|}\,\d \sigma(\y)\\
     &= \int_{S_1(\vO)}   \frac{  u(\x) }{|S_1(\vO)|}\,\d \sigma(\y) = v(\x,0)
\end{aligned}
\end{equation*}
i.e. $t\mapsto v(\x,t)$ is continuous from the right at $t=0$. Moreover, by~(\ref{repv02})
\begin{equation*}
\begin{aligned}
   v(\x,\tau+)  
     &= \lim_{s\to 0,\,s>0}\int_{S_1(\vO)}   \frac{  v(\x +R(s)\,\y,\tau) }{|S_1(\vO)|}\,\d \sigma(\y)\\
     &= \int_{S_1(\vO)}   \frac{  v(\x,\tau) }{|S_1(\vO)|}\,\d \sigma(\y) = v(\x,\tau)\,.
\end{aligned}
\end{equation*} 
In summary, we have shown that $t\mapsto v(\x,t)$ is continuous for $t\in (0,\tau]$.  \\
Assume that $t\mapsto v(\x,t)$ is continuous for $t\in (0,\tau_m]$. \\ 
Let $t\in (\tau_m,\tau_{m+1})$. Analogously as above, it follows by induction from 
\begin{equation*}
\begin{aligned}
   &| v(\x,t+s)- v(\x,t)|    \\
     &\qquad \leq   \int_{S_1(\vO)}  \frac{| v(\x+R(t+s)\,\y,\tau_m)
                       - v(\x+R(t)\,\y,\tau_m)| }{|S_1(\vO)|} \,\d \sigma(\y)
\end{aligned}
\end{equation*}
and~(\ref{repv02}) that $t\mapsto v(\x,t)$ is continuous for $t\in (\tau_m,\tau_{m+1}]$. 
This proves the claim of part a). \\ 
b): Let $u$ be differentiable. From~(\ref{repv02}), it follows  
\begin{equation*}
\begin{aligned}
   \frac{\partial v}{\partial t} (\x,t)
      &= \frac{\int_{S_1(\vO)} [\nabla v(\x+R(t)\,\y,\tau_{n(t)})] \cdot \y\,\d \sigma(\y) }
                           {|S_1(\vO)|}\,c_0\,,
\end{aligned}
\end{equation*}
i.e.  $\frac{\partial v}{\partial t} (\x,\cdot)$ is continuous on $(\tau_{n(t)},\tau_{n(t)+1})$ 
and 
\begin{equation}\label{vt+}
\begin{aligned}
      \frac{\partial v}{\partial t}(\cdot,\tau_m+) = 0     \qquad (m\in\N) \,,
\end{aligned}
\end{equation}
since $R(\tau_m+)=\lim_{s\to 0+} R(\tau_m+s)=0$ and $\int_{S_1(\vO)} \y\,\d \sigma(\y)=\vO$. 
Moreover, Proposition~\ref{prop:green02} implies that
\begin{equation}\label{vt-}
\begin{aligned}
   \frac{\partial v}{\partial t}(\cdot,\tau_{m+1}-)  
         = u *_\x \frac{\partial G}{\partial t}(\cdot,\tau_{m+1}-) 
\qquad (m\in\N)\,
\end{aligned}
\end{equation}
and therefore~(\ref{jumpcond}) holds. This concludes the proof.
\end{proof}

\begin{figure}[!ht]
\begin{center}
\includegraphics[height=8.0cm,angle=0]{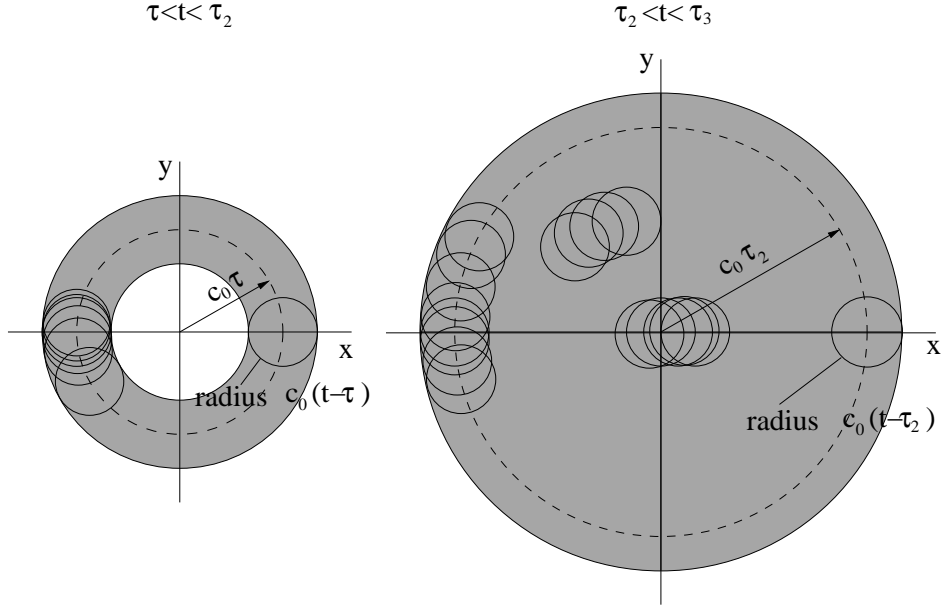}
\end{center}
\caption{Left figure: If $N=2$ and $t\in (\tau,\tau_2)$, then $\supp(G(\cdot,t))$ 
is an annulus generated by circles of radius $c_0\,(t-\tau)$ and center at $\x'\in S_{c_0\,\tau}(\vO)$. 
Right figure:  If $N=2$ and $t\in (\tau_2,\tau_3)$, then $\supp(G(\cdot,t))$ 
is a disc generated by circles of radius $c_0\,(t-\tau_2)$ and center at $\x'\in B_{c_0\,\tau_2}(\vO)$.
}
\label{fig:SuppG2}
\end{figure}

\begin{prop}\label{prop:vcont2}
Let $v$ and $v_c$ be defined as in Definition~\ref{def:diff2} with $u_c\in L^1(\R^N)$ for each 
$c\in (0,c^*]$. Moreover, let  $\x\in\R^N$ and $t>0$.\\
a) If $(\x,c)\mapsto u_c(\x)$ is continuous on $\R^N\times (0,c^*)$, then $c\mapsto v_c(\x,t)$ is 
continuous on $(0,c^*)$. \\
b) If $c\mapsto f(c)$ and $c\mapsto u_c(\x)$ are continuous on $(0,c^*)$, then  
$t\mapsto \frac{\partial v}{\partial t} (\x,t)$ is continuous on $(\tau^*,\infty)$. 
\end{prop}

\begin{proof}
a) Let $c\in (0,c^*)$ and $\tau:=\lambda/c$.  We have to show that
\begin{equation}\label{limvc}
   \lim_{\epsilon\to 0+}  v_{c-\epsilon}(\x,t) =  \lim_{\epsilon\to 0+} v_{c+\epsilon}(\x,t)
\end{equation}
for the three cases 
$$
    t<\tau\,,  \qquad \qquad \quad  t=\tau   \qquad\qquad \mbox{and} \qquad\qquad  t>\tau \,.
$$
We prove only the case $t=\tau$. The remaining two cases are proven similarly. 
Let $c_\epsilon:=c-\epsilon$ and $\tau_\epsilon:=\lambda/c_\epsilon$. 
From~(\ref{repv01}) with $t=\tau<\tau_\epsilon$ and the continuity of $(\x,c)\mapsto u_c(\x)$,  
it follows that 
\begin{equation*}
\begin{aligned}
  \lim_{\epsilon\to 0+}  v_{c_\epsilon}(\x,t) 
     &= \lim_{\epsilon\to 0+} \int_{S_1(\vO)}  \frac{v_{c_\epsilon}(\x+c_\epsilon\,t,0)}{|S_1(\vO)|} \,\d \sigma (\y)\\
     &= \int_{S_1(\vO)}  \frac{u_c(\x+c\,t)}{|S_1(\vO)|} \,\d \sigma (\y) 
      = v_c(\x,t)\,.
\end{aligned}
\end{equation*}
Let $c_\epsilon:=c+\epsilon$ and $\tau_\epsilon:=\lambda/c_\epsilon$. 
From~(\ref{repv01}) and~(\ref{repv02}) with $t=\tau>\tau_\epsilon$ and the continuity of $(\x,c)\mapsto u_c(\x)$,  
it follows that 
\begin{equation*}
\begin{aligned}
  \lim_{\epsilon\to 0+}  v_{c_\epsilon}(\x,t) 
     &= \lim_{\epsilon\to 0+} \int_{S_1(\vO)}  \frac{v_{c_\epsilon}(\x+c_\epsilon\,(t-\tau_\epsilon),\tau_\epsilon)}
                   {|S_1(\vO)|}\,\d \sigma (\y)\\
     &= \lim_{\epsilon\to 0+} \int_{S_1(\vO)}\int_{S_1(\vO)} 
              \frac{u_{c_\epsilon}(\x + c_\epsilon\,(t-\tau_\epsilon)\,\y_1 + c_\epsilon\,\tau_\epsilon\,\y_2)}
                                   {|S_1(\vO)|^2} \,\d \sigma (\y_1)\,\d \sigma (\y_2) \\
     &= \int_{S_1(\vO)}\int_{S_1(\vO)} 
              \frac{u_c(\x + 0 + c\,\tau\,\y_2)}
                                   {|S_1(\vO)|^2} \,\d \sigma (\y_1)\,\d \sigma (\y_2) 
     = v_c(\x,\tau) 
     = v_c(\x,t)\,,
\end{aligned}
\end{equation*}
which proves~(\ref{limvc}) for $t=\tau$. \\
b) According to Definition~\ref{def:diff2}, we have 
\begin{equation*}
\begin{aligned}
    v(\x,t) =\int_{\tau^*}^\infty w(\x,t,\tau)  \,\d \tau 
\qquad\mbox{with}\qquad
      w(\x,t,\tau) =\frac{\lambda}{\tau^2}\,f\left(\frac{\lambda}{\tau}\right)\,v_{\lambda/\tau}(\x,t)\,.
\end{aligned}
\end{equation*}
For $t\in (\tau^*,\infty)$, let $ \tau_0:=\tau^*$ and $\tau_m:=m\,t$ for $m\in\N$. 
From~(\ref{limvc}),  we get  
\begin{equation*}
\begin{aligned}
     \left[w(\cdot,t,s) \right]_{s=\tau_m-}^{s=\tau_m+}
          = \frac{\lambda}{t^2}\,f\left(\frac{\lambda}{t}\right)\, 
                 \left[v_{\lambda/s}(\cdot,t)\right]_{s=\tau_m-}^{s=\tau_m+} 
          = 0   \qquad\quad (m\in\N)\,,
\end{aligned}
\end{equation*}
and thus 
\begin{equation*}
\begin{aligned}
    \frac{\partial v}{\partial t}(\x,t) 
        &= \sum_{m=0}^\infty \int_{\tau_m}^{\tau_{m+1}} \frac{\partial w}{\partial t}(\x,t,\tau)  \,\d \tau 
           - \sum_{m=1}^\infty  \left[w(\x,t,s) \right]_{s=\tau_m-}^{s=\tau_m+}  \\
        &= \sum_{m=0}^\infty \int_{\tau_m}^{\tau_{m+1}} \frac{\partial w}{\partial t}(\x,t,\tau)  \,\d \tau\,,
\end{aligned}
\end{equation*}
where $w(\x,t,\tau)$ is differentiable for each $\tau\in (\tau_m,\tau_{m+1})$. 
This concludes the proof. 
\end{proof}

In the following we focus on the case $f(c)=\delta(c-c_0)$ and usually write $G$ instead of
$G_{c_0}$. We now specify the support of the Green function $G$ and show that
causality condition~(\ref{causaGcond}) is satisfied.

\begin{prop}\label{prop:causality}
Let $\tau$, $c$ be as in A1)-A2) and $f(c)=\delta(c-c_0)$. Moreover, let $G$ denote the
Green function of diffusion, i.e. $G$ satisfies~(\ref{defcausdiff1}) with $u=\delta(\x)$.\\
a) Let $N\in\N$. Then $G$ satisfies causality condition~(\ref{causaGcond}) and
\begin{equation}\label{suppG1}
\begin{aligned}
   \supp\,(G(\cdot,t)) = S_{c_0\,t}(\vO)
   \qquad \mbox{for}\qquad t\in [0,\tau]\,.
\end{aligned}
\end{equation}
b) If $N>1$, then
\begin{equation}\label{suppG2}
\begin{aligned}
   \supp\,(G(\cdot,t)) = \overline{\B_{c_0\,t}(\vO) \backslash \B_{c_0\,(2\,\tau-t)}(\vO)}
   \qquad \mbox{for}\qquad t\in (\tau,2\,\tau]\,
\end{aligned}
\end{equation}
and
\begin{equation*}
\begin{aligned}
   \supp\,(G(\cdot,t)) = \overline{\B_{c_0\,t}(\vO)}
   \qquad \mbox{for}\qquad t \geq 2\,\tau\,.
\end{aligned}
\end{equation*}
\end{prop}

\begin{proof}
a) Identity~(\ref{suppG1}) follows at once from~(\ref{repG1}). 
To prove the rest of the claim, we perform a proof by induction. 
According to~(\ref{suppG1}) we have for $t\in [0,\tau]$:
\begin{equation}\label{helpGB}
    \supp(G(\cdot,t))\subseteq \overline{\B_{c_0\,t}(\vO)} \,. 
\end{equation}
Assume that~(\ref{helpGB}) holds for $t\in [0,\tau_m]$ (Induction assumption). 
From~(\ref{mainpropG}) with $t\in (\tau_m,\tau_{m+1}]$ and~(\ref{suppG1}) together with the support 
theorem of distributions with compact support (cf. Theorem~4.3.3 in~\cite{Hoe03}), it follows 
for $t\in (\tau_m,\tau_{m+1}]$ that 
\begin{equation*}
\begin{aligned}
      \supp (G(\cdot,t)) 
           &= ch\,\supp (G(\cdot,\tau_m)) + ch\,\supp (G(\cdot,t-\tau_m)) \\
           &\subseteq \overline{\B_{c_0\,\tau_m}(\vO)} + \overline{\B_{c_0\,(t-\tau_m)}(\vO)} 
             = \overline{\B_{c_0\,t}(\vO)}  \,,
\end{aligned}
\end{equation*}
i.e.~(\ref{helpGB}) holds. By induction, we infer that~(\ref{helpGB}) 
holds for $t\geq 0$, i.e. causality condition~(\ref{causaGcond}) is satisfied. \\
b)
To prove~(\ref{suppG2}) we note that~(\ref{mainpropG}) and~(\ref{repG1}) imply for 
$t\in(\tau,2\,\tau]$:
\begin{equation*}
\begin{aligned}
   G(\x,t)
     = \frac{  \int_{\R^N} \delta(c_0\,\tau - |\x'|) \, \delta(c_0\,(t-\tau) - |\x-\x'|) \,\d \x'
             }{|S_{c_0\,\tau}(\vO)|\,|S_{c_0\,(t-\tau)}(\vO)|}\,.
\end{aligned}
\end{equation*}
Thus, for this case, the support of $G(\cdot,t)$ is the union of the set of spheres with 
radius $c_0\,(t-\tau)$ and center $\x'\in S_{c_0\,\tau}(\vO)$ (cf. Fig.~\ref{fig:SuppG2} for $N=2$), 
which is the shell $\overline{\B_{c_0\,t}(\vO) \backslash \B_{c_0\,(2\,\tau-t)}(\vO)}$.

If $t\in (\tau_2,\tau_3]$, then according to~(\ref{mainpropG}) and~(\ref{suppG2}) $\supp(G(\cdot,t))$ 
is the union of the set of spheres with radius $c_0\,(t-\tau_2)$ and center 
$\x'\in \overline{B_{c_0\tau_2}(\vO)}$ (cf. Fig.~\ref{fig:SuppG2} for $N=2$), 
which is the closed ball $\overline{\B_{c_0\,t}(\vO)}$. 
By induction, it follows from~(\ref{mainpropG}) for $t\in (\tau_m,\tau_{m+1}]$ with $m\geq 2$
that $\supp(G(\cdot,t))$ is the union of the set of spheres with radius $c_0\,(t-\tau_m)$ 
and center $\x'\in \overline{B_{c_0\tau_m}(\vO)}$, which is the closed ball $\overline{\B_{c_0\,t}(\vO)}$. 
This concludes the proof. 
\end{proof}

\subsection{Interpretation of the  diffusion model}
\label{subsec-interpr}

We now give two interpretations of diffusion model~(\ref{defcausdiff1}). Firstly, we interpret 
it as transport process and secondly we interpret it as stochastic process.

\subsubsection*{Interpretation as transport process}

Let $G$ denote the Green function of the diffusion process~(\ref{defcausdiff1}) that takes places 
with constant speed $c_0$. 
Then the concentration of the diffusing substance is given by 
$$
   v(\x,t) = (G(\cdot,t) *_\x u)(\x)  \qquad  \qquad v(\x,0)=u(\x)\,
$$
and thus it is sufficient to explain the process for the case $u(\x)=\delta(\x)$. 
Let $u(\x)=\delta(\x)$. Then the concentration in point $\x=\vO$ splits up and starts to 
spread out in each direction at time $t=0$ and propagates with speed $c_0$ during the time 
period $(0,\tau)$. 
Because of Lemma~\ref{lemm:masscons}, mass is conserved during this transport process.
According to relation~(\ref{mainpropG}), we have
$$
    G(\x,t) =  \int_{\R^N} G(\x,\tau)\,G(\x-\x',t-\tau)\,\d \x'
       \qquad \mbox{for}\qquad t\in(\tau,\tau_2)\,,
$$
and hence the concentration on the sphere $S_{c_0\,\tau}(\vO)=\supp\,(G(\x,\tau))$ 
(cf.~(\ref{repG1})) splits up and starts to spread out in each direction at time $t=\tau$. 
The concentration propagates with speed $c_0$ during the time period $(\tau,\tau_2)$ 
(cf. Fig.~\ref{fig:SuppG2} and Proposition~\ref{prop:causality}) 
such that the total mass is conserved. 
This process of splitting up and spreading out carries on forever if the medium is unbounded.

We note that this model allows particles to cross (without interfering) and 
consequently not each space-time point can be associated a single velocity vector. 
Indeed, from~(\ref{defj}) we see that the flux density $\vj(\x,\tau_{m+1})$ contains the 
velocity vectors  
$$
                 \vc = c_0\,\y  \qquad\mbox{for}\qquad \y\in S_1(\vO)  
   \qquad\mbox{if}\qquad   v(\x+c_0\,\tau\,\y,\tau_m) \not= 0\,.
$$

\subsubsection*{Interpretation as stochastic process}

We can give a probabilistic interpretation of the diffusion model, too. 
We interpret $v(\cdot,t)$ as the probability density of a single particle at time $t$ 
which is part of a diffusion process satisfying~(\ref{defcausdiff1}). Mathematically, 
this means that 
$$
    P_t(\Omega) := \int_\Omega v(\x',t)\,\d \x'
$$
describes the probability that the particle lies in the region $\Omega$ at time $t$. 
For simplicity, we assume $N=3$. Let $t\in (\tau_m,\tau_{m+1}]$ and 
\begin{itemize}
\item [1)] $A_t$ denote the event that this particle arrives (within a fixed cube $\Omega_\x$ 
           with side length $\d l$) around point $\x$ at time $t$ and

\item [2)] $B_m$ denote the event that this particle arrives (within a fixed spherical shell 
           $\Gamma_{\x,m}$ with thickness $\d R$)\footnote{The length $\d l$ of the cube determines the thickness 
           $\d R$ of the shell by $\d l^3 = |S_{\d R}(\x)| \,\d R$ (origin is at $\x$).} 
           around surface $S_{R(t)}(\x)$ at time $\tau_m$.
\end{itemize}
If event $A_t$ occurs, then the probability of the particle 
to be outside the spherical shell $\Gamma_{\x,m}$ at time $\tau_m$ is zero by~(\ref{defcausdiff1}). 
Moreover, if event $A_t$ occurs and the particle was around point $\x'\in S_{R(t)}(\x)$ at time 
$\tau_m$, then the particle was \emph{around}
$$
   \y(s) = \x' + (\x-\x')\,\frac{s-\tau_m}{t-\tau_m}
        \qquad \mbox{for time $s\in (\tau_m,\tau_{m+1})$}\,.
$$
However, we do not exactly know where the particle was during this period.
In Proposition~\ref{prop:causality} we showed that this type of probabilistic process
guarantees causality and thus in contrast to the noncausal models of stochastic processes, 
the particle lies within a \emph{bounded region} for each time point.

\subsection{The equation of diffusion}

Although \emph{continuity equation}~(\ref{gconteq}) with flux density~(\ref{defj}) can be 
considered as \emph{the} equation of diffusion, we derive an alternative equation that is a proper  
partial differential equation. We refer only to this equation as \emph{the equation of diffusion}.

According to Proposition~\ref{prop:green02} and Theorem~\ref{theo:propG}, the  
concentration of a substance diffusing with constant speed, say $c_0$, is given by
$$
    v(\cdot,t) = u *_\x G(\cdot,\tau_{n(t)})*_\x G (\cdot,t-\tau_{n(t)})  
     \qquad\quad t\geq 0\,,
$$
where $u$ and $G$ denote the initial concentration and the Green function, respectively. 
In order to derive the equation of diffusion, we utilize the following identities for 
$t\in (\tau_m,\tau_{m+1}]$ with $m\in\N_0$:
\begin{equation}\label{helpwaveeq}
\begin{aligned}
    &\nabla^2 v(\cdot,t)
        = u *_\x G(\cdot,\tau_m)*_\x  \nabla^2 G (\cdot,t-\tau_m)  \\
    &\frac{\partial^k v}{\partial t^k} (\cdot,t)
        = u *_\x G(\cdot,\tau_m)*_\x \frac{\partial^k G}{\partial t^k} (\cdot,t-\tau_m)\,,
\end{aligned}
\end{equation}
where $k\in\{1,\,2\}$. 
We see if $G|_{\R^N\times (0,\tau)}$ satisfies a partial differential equation with time dependend 
coefficients $a(t)$, $b(t)$,$\,\ldots\,$, then $v$ satisfies the same equation with time dependend
coefficients $a(t-\tau_{n(t)})$, $b(t-\tau_{n(t)})$,$\,\ldots$ as long as $t\not=\tau_{n(t)}$. 
In addition, we have to pay regard that $\frac{\partial v}{\partial t}$ has at $t=\tau_m$ a jump 
of size  (cf. Proposition~\ref{prop:vcont})
\begin{equation}\label{helpwaveeqb}
   \left[\frac{\partial v}{\partial t}\right]_{t=\tau_m-}^{t=\tau_m+}
       = -u *_\x \frac{\partial G}{\partial t}(\cdot,\tau_m-)\,.
\end{equation}

For the derivation of the equation of diffusion we need two lemmata. 

\begin{lemm}\label{lem:propdelta2}
Let $R(t)$ be as in A2), then
\begin{equation*}\label{propdelta2}
\begin{aligned}
  \frac{\delta'(R(t)-|\x|)}{|\x|}
      = \frac{\delta'(R(t)-|\x|)}{R(t)}
       - \frac{\delta(R(t)-|\x|)}{R(t)^2}\,.
\end{aligned}
\end{equation*}
\end{lemm}

\begin{proof}
Without loss of generality we assume that $c_0=1$, i.e. $R(t)=t$. 
Let $C_0^\infty(\R^{N+1})$ denote the space of $C^\infty-$function with compact support in $\R^{N+1}$. 
For $\psi\in C_0^\infty(\R^{N+1})$, we have to show that
\begin{equation*}
\begin{aligned}
    \int_{\R^{N+1}} \psi(\x,t)\, &(R(t)-|\x|)\, \delta'(R(t)-|\x|)\, \d \x \,\d t \\
         &=  -\int_{\R^{N+1}} \psi(\x,t)\, \frac{|\x|}{R(t)}\, \delta(R(t)-|\x|)\, \d \x \,\d t \,.
\end{aligned}
\end{equation*}
But this is equal to 
\begin{equation*}
\begin{aligned}
    &-\int_{\R^{N+1}} \left[\frac{\partial\psi}{\partial t}(\x,t)\, (R(t)-|\x|)\, 
                         + \psi(\x,t)  \right]_{R(t)=|\x|}\, \d \x  \\
         &=  -\int_{\R^{N+1}} \psi(\x,|\x|)\,\, \d \x \,
\end{aligned}
\end{equation*}
which holds for every $\psi\in C_0^\infty(\R^{N+1})$. 
\end{proof}

\begin{lemm}\label{lem:waveqg}
Let $\tau$, $c$, $\tau_m$, $R(t)$ be as in A1)-A2) and let $G$ be defined in~(\ref{defcausdiff1}) 
with $u=\delta(\x)$. Then $g(\cdot,t):=G(\cdot,t-\tau_m)$ ($m\in\N_0$) solves the following  
equation on $\R^N\times (\tau_m,\tau_{m+1})$
\begin{equation}\label{seqofeq0}
\begin{aligned}
   &\frac{1}{c_0^2}\frac{\partial^2 g}{\partial t^2}
         +  \frac{1}{c_0}\,\frac{(N-1)}{R(t)}\,\frac{\partial g}{\partial t} 
    - \nabla^2 g
         = 0   
\end{aligned}
\end{equation}
with initial conditions 
\begin{equation}\label{seqofeq0b}
\begin{aligned}
   &g(\x,\tau_m) = \delta(\x)\,\qquad\mbox{and}\qquad
       \frac{\partial g}{\partial t}(\x,\tau_m) = 0   \qquad \x\in\R^N\,.
\end{aligned}
\end{equation}
\end{lemm}

\begin{proof}
We note that $g(\x,t)=G(\x,R(t)/c_0)$ for $\x\in\R^N$ and $t\in (\tau_m,\tau_{m+1})$. 
From~(\ref{repG1}),~(\ref{propS}) and Lemma~\ref{lem:propdelta2}, we get
\begin{equation*}
\begin{aligned}
     \nabla^2 g (\x,t)
         &= \frac{\delta''(R(t)-|\x|)}{|S_{R(t)}(\vO)|}
           - \frac{N-1}{|\x|}\,\frac{\delta'(R(t)-|\x|)}{|S_{R(t)}(\vO)|} \\
         &=  \frac{\delta''(R(t)-|\x|)}{|S_{R(t)}(\vO)|}
           - (N-1)\,\left[ \frac{\delta'(R(t)-|\x|)}{|S_{R(t)}(\vO)| \,R(t)}
               - \frac{g(\x,t)}{R(t)^2} \right]\,.
\end{aligned}
\end{equation*}
Moreover,~(\ref{repG1}) and~(\ref{propS}) imply 
\begin{equation*}
\begin{aligned}
     \frac{\partial g}{\partial t}(\x,t)
         = c_0\,\frac{\delta'(R(t)-|\x|)}{|S_{R(t)}(\vO)|}
           - \frac{c_0\,(N-1)}{R(t)}\,g(\x,t) \,
\end{aligned}
\end{equation*} 
and
\begin{equation*}
\begin{aligned}
     \frac{\partial^2 g}{\partial t^2}(\x,t)
         = c_0^2\,\frac{\delta''(R(t)-|\x|)}{|S_{R(t)}(\vO)|}
           - \frac{2\,c_0^2\,(N-1)}{R(t)}\,\frac{\delta'(R(t)-|\x|)}{|S_{R(t)}(\vO)|}
           + \frac{c_0^2\,(N-1)\,N}{R(t)^2}\,g(\x,t) \,.
\end{aligned}
\end{equation*}
From these identities, we obtain equation~(\ref{seqofeq0}) for $t\in (\tau_m,\tau_{m+1})$\,.

Since $G(\x,0) = \delta(\x)$ by~(\ref{repG1b}), we get 
$$
       g(\x,\tau_m) = G(\x,0) = \delta(\x)\,,
$$
which proves the first initial conditions in~(\ref{seqofeq0b}).  
The second initial conditions in~(\ref{seqofeq0b}) follow from Proposition~\ref{prop:vcont}.
\end{proof}

To complete the derivation of the diffusion equation we have to incorporate the 
jump condition~(\ref{helpwaveeqb}) via a source term (right hand side term).

\begin{prop}\label{prop:causalwaveeq}
Let $c_0>0$ be constant, $R(t)$ be as in A2) and $v$ be defined as in~(\ref{defcausdiff1})
with differentiable $u\in L^1(\R^N)$. Then $v$ satisfies equation\footnote{To be precise, 
we use the generalized time derivative, denoted by $\partial_t$, since 
$v$ is not differentiable at $t=\tau_m$.}
\begin{equation}\label{seqofeq1}
\begin{aligned}
   &\frac{1}{c_0^2}\,\partial_t^2 v 
         + \frac{1}{c_0}\,\frac{(N-1)}{R(t)}\,\partial_t v 
    - \nabla^2 v = f[u] \,\quad \mbox{on $\R^N\times (0,\infty)$}
\end{aligned}
\end{equation}
with source term
\begin{equation}\label{seqofeq1a}
\begin{aligned}
   f[u](\cdot,t) = -\frac{1}{c_0^2}\,u *_\x 
           \sum_{m\in\N} \frac{\partial G}{\partial t}(\cdot,\tau_m-)\,\delta(t-\tau_m)
\end{aligned}
\end{equation}
and initial conditions 
\begin{equation}\label{seqofeq1b}
\begin{aligned}
   &v(\cdot,0) = u\,\qquad\mbox{and}\qquad
       \frac{\partial v}{\partial t}(\cdot,0) = 0 \,.
\end{aligned}
\end{equation}
\end{prop}

\begin{proof}
All except the source term $f$ follows at once from Lemma~\ref{lem:waveqg} with~(\ref{helpwaveeq}).  
According to~(\ref{helpwaveeqb}) (cf. Proposition~\ref{prop:vcont}), $\frac{\partial v}{\partial t}$ 
has a jump at $t=\tau_m$ and thus 
\begin{equation*}
\begin{aligned}
    \partial_t^2 v 
     &=  \frac{\partial^2 v}{\partial t^2} 
         + \sum_{m\in\N}\left[\frac{\partial v}{\partial t}\right]_{t=\tau_m-}^{t=\tau_m+}\,
                          \delta(t-\tau_m)  \\
     &=  \frac{\partial^2 v}{\partial t^2} 
       -u *_\x \sum_{m\in\N} \frac{\partial G}{\partial t}(\cdot,\tau_m-)\,\delta(t-\tau_m)\,,
\end{aligned}
\end{equation*}
where $\partial_t$ and $\frac{\partial}{\partial t} $ denote the generalized and the pointwise 
time derivative, respectively. From this and~(\ref{seqofeq0}), we infer equation~(\ref{seqofeq1}) 
with the source term~(\ref{seqofeq1a}). 
\end{proof}

\begin{rema}\label{rema:limit}
The classical diffusion equation 
\begin{equation*}
\begin{aligned}
    \frac{\partial v}{\partial t} 
      - D_0\,\nabla^2 v  = 0
    \qquad\mbox{with}\qquad v(\x,0) = u(\x)
\end{aligned}
\end{equation*}
can be obtained via a special limit from model~(\ref{defcausdiff1}) for the time points 
$t=\tau_m$ for $m\in\N_0$ with $\tau\to 0$ under the side condition 
$$
            D_0:=\frac{c_0^2\,\tau}{2\,N}=const\,.
$$
Note that $D_0 = D(\tau/2) = \frac{1}{\tau}\,\int_0^\tau D(s)\,\d s$ with $D(s)$ defined as in~(\ref{Dt}).  
For these time points, representation~(\ref{defcausdiff1}) can be written as follows 
\begin{equation*}
\begin{aligned}
  \int_{S_1(\vO)} v(\x,\tau_{m+1}) \, \d \sigma(\y)
      = \frac{1}{2}\int_{S_1(\vO)} [v(\x+R(\tau)\,\y,\tau_m) 
                             + v(\x-R(\tau)\,\y,\tau_m)] \,\d \sigma(\y)
\end{aligned}
\end{equation*}
which is equivalent to
\begin{equation*}
\begin{aligned}
  \int_{S_1(\vO)} & \frac{v(\x,\tau_{m+1}) - v(\x,\tau_m) }{\tau} \, \d \sigma(\y)  \\
      &= N\,D_0\,\int_{S_1(\vO)} \frac{v(\x+R(\tau)\,\y,\tau_m) 
                  - 2\,v(\x,\tau_m) + v(\x-R(\tau)\,\y,\tau_m)}{R(\tau)^2} \,\d \sigma(\y)\,.
\end{aligned}
\end{equation*}
If we perform the limit $\tau\to 0$ such that $D_0=const.$ (and consequently $c_0\to\infty$), 
then we obtain from~(\ref{Dkl}) 
\begin{equation*}
\begin{aligned}
   \,\frac{\partial v}{\partial t}(\x,t) 
      = \frac{N\,D_0}{|S_1(\vO)|}\,\int_{S_1(\vO)} \lim_{s\to 0}
            \frac{\partial^2 v(\x+s\,\y,t)}{\partial s^2} \,\d \sigma(\y) 
      = \sum_{k,l=1}^3 D_{k,l}\,\frac{\partial^2 v(\x,t)}{\partial x_k\partial x_l} 
\end{aligned}
\end{equation*}
with $D_{k,l}= D_0\,\,\delta_{k,l}$ for $k,\,l\in\{1,2,\ldots,N\}$.  
Apart from the causality problem, this limit is not reasonable, since 
$\frac{\partial v}{\partial t}$ is not continuous at $t=\tau_m$ (cf. Proposition~\ref{prop:vcont}). 
\end{rema}

\section{Conclusions}
\label{sec-conclusion}

Since it is not possible that a \emph{causal} diffusion process satisfies the (strongly 
continuous) semigroup property, we developed a causal model of diffusion with constant 
speed $c$ which satisfies the semigroup property only at a discrete set $M_c$ of time points. 
In contrast to the classical and noncausal diffusion model, the causal diffusion model 
is not differentiable at the discrete set $M_c$ of time points.
This property of non-smoothness forbids the limit process that transforms the 
discrete and causal diffusion model from Physics and Stochastics into the continuous and 
noncausal standard diffusion model. 
It is surprising that a smoothness assumption, which is too strong, has an effect on the 
causality of a model.  
Another consequence of causality is that diffusion (with constant speed) satisfies an 
inhomogeneous wave equation with a time dependent coefficient. 
We have seen that a diffusion process with variing speed $c\in (0,c^*]$ is a superposition 
of processes each of which satisfying the semigroup property on a 
discrete set $M_c$ such that $\bigcup_{c\in (0,c^*]} M_c= [\tau^*,\infty)$. However, the 
total process does not satisfy the strongly continuous semigroup property on $[\tau^*,\infty)$ 
and is not $C^\infty$ with respect to time.

Finally, I would like to give two remarks. \\
1) It is a valid option to model diffusion-like or stochastic processes by vector-valued 
semigroups. 
An interesting example from physics may be the Dirac equation (cf.~\cite{Wei95,Zee03}). 
Whether such an approach is more convenient or powerful, then modeling with 
hyperbolic equations or ``other equations'' is not clear.  \\
2) From the mathematical point of view, it is interesting whether a causal variant of the 
Schr\"odinger equation is possible such that for example quantum tunneling obeys causality.  
It seems natural to the author that the unresolved causality problem of diffusion had an  
influence on the development of quantum mechanics and consequently some of the causality 
problems in this theory may be resolvable.  
Moreover, it is intriguing that the existence of antiparticles in quantum field theory is 
based on a ``special form'' of causality (cf. Chapter 2 Section 13 in~\cite{Wei72}).

\section{Appendix}
\label{sec-app}

In this appendix we clarify our notion of Fourier transform, convolution, formulate the 
Paley-Wiener Schwartz Theorem which is used for the causality analysis, and 
complete the proof of Theorem~\ref{theo:main} from Section~\ref{sec-semigroups}.

\subsection{Definition of the Fourier transform}

The Fourier transform of $f\in L^1(\R^N)$ is defined by
\begin{equation*}
\begin{aligned}
   \F\{f\}(\k) := (2\,\pi)^{-N/2}\,\int_{\R^N} e^{\i\,\k\cdot \x} \,f(\x)\,\d \x
  \qquad\quad \k\in\R^N\,,
\end{aligned}
\end{equation*}
where $\k\cdot \x:=\sum_{j=1}^N k_j\,x_j$. In this notion the convolution theorem reads 
as follows
\begin{equation*}
\begin{aligned}
    \F\{f\}\,\F\{g\} = (2\,\pi)^{-N/2} \,\F\{f *_\x g\} 
\qquad f,\,g\in L^1(\R^N), 
\end{aligned}
\end{equation*}
where the space-convolution is defined by
$$
   (f *_\x g)(\x) := \int_{\R^N} f(\x')\,g(\x-\x')\,\d \x'\,. 
$$

\subsection{The Paley-Wiener Schwartz Theorem}

\vspace{0.5cm}
We note that the \emph{supporting function} of a compact set $K\in\R^N$ is defined by
\begin{equation}\label{supfunc}
\begin{aligned}
   H_K(\xi) := \sup_{\x\in K} \,\sum_{j=1}^N x_j\,\xi_j
     \qquad\mbox{for}\qquad \xi\in\R^N\,.
\end{aligned}
\end{equation}
Moreover, 
$$
    \mbox{Im}(z):=(\mbox{Im}(z_1),\,\ldots,\,\mbox{Im}(z_N)) 
\qquad\mbox{and}\qquad 
    |z|^2:=\sum_{j=1}^N z_j^2 \qquad z\in\C^N\,.
$$
By $B_R(\vO)$ we denote the ball of radius $R$ and center $\vO$. 
The Paley-Wiener Schwartz Theorem (cf. Theorem~7.3.1 in~\cite{Hoe03}) reads as follows

\begin{theo}\label{th:PWS1}
Let $N\in\N$. $f:\R^N\to\R$ is a distribution with support in the compact set 
$K\subseteq\R^N$ if and only if
\begin{itemize}
\item [P1)] $\k\in\R^N\mapsto \F\{f\}(\k)\in\R$ can be extended to an entire function 
            denoted by $z\in\C^N\mapsto\F\{f\}(z)\in\C$ and

\item [P2)] there exist $M>0$ and $C>0$ such that
$$
         |\F\{f\}(z)| \leq C\,(1+|z|)^M\,e^{H_K(\mbox{Im}(z))}
  \qquad\quad\mbox{for all $z\in\C^N$}\,.
$$
\end{itemize}
The extension $\F\{f\}(z)$ is the Fourier-Laplace transform of $f$, and $M$
can be chosen as the order of the distribution $f$.
\end{theo}

\begin{rema}\label{rema:PWS}
a) If $K= [-c_0,c_0]^N$, then
\begin{equation*}
\begin{aligned}
   H_K(\mbox{Im}(z)) = c_0\,\sum_{j=1}^N |\mbox{Im}(z_j)|
     \qquad\mbox{for}\qquad z\in\C^N\,
\end{aligned}
\end{equation*}
and if $K=B_{c_0}(\vO)$, then
\begin{equation*}
\begin{aligned}
   H_K(\mbox{Im}(z)) = c_0\,\sqrt{\sum_{j=1}^N \mbox{Im}(z_j)^2}
     \qquad\mbox{for}\qquad z\in\C^N\,.
\end{aligned}
\end{equation*}
b) Let $K\subseteq B_{c_0}(\vO)$ and 
\begin{equation}\label{defwz}
\begin{aligned}
   w(z_1):=(z_1,0,\ldots,0)\in\C^N \,.
\end{aligned}
\end{equation}
Then from property P2) in Theorem~\ref{th:PWS1}, it follows that
\begin{equation*}
\begin{aligned}
      \limsup_{r\to\infty} \frac{\log \max_{|z_1|=r}|\F\{f\}(w(z_1))|}{r} \leq c_0\,,
\end{aligned}
\end{equation*}
i.e. $\F\{f\}\circ w$ is a function of exponential type $\leq c_0$ (cf. Definition 2.1.3
and Theorem 6.8.1 in~\cite{Boa73}). \\
c) Let $f(\x)$ be rotational symmetric in $\x\in\R^N$ and $\F\{f\}(z)$ be entire. 
Then  
$\F\{f\}(z)$ is rotational symmetric in $z\in\C^N$, too, and property P2) of 
Theorem~\ref{th:PWS1} holds if and only if the following property holds:
\begin{itemize}
\item [P2')] there exist $M>0$ and $C>0$ such that
$$
         |\F\{f\}(w(z_1))| \leq C\,(1+|z_1|)^M\,e^{H_K(\mbox{Im}(z_1))}
  \qquad\quad\mbox{for all $z_1\in\C$}\,.
$$
\end{itemize}
Here $w$ is defined as in b). 
\end{rema}

\subsection{Final parts of the proof of Theorem~\ref{theo:main}}

Let $\hat\a(z)$ be entire, $G$ be defined as in~(\ref{defG}) and 
$$
     f_1(z) := G(\cdot,1) = \exp\{-\hat a(z)\}\,. 
$$
To complete the proof of Theorem~\ref{theo:main} in Section~\ref{sec-semigroups}, we have 
to prove two statements:
\begin{itemize}
\item [A)] If $\hat\a$ is a polynomial in $\i\,z$ of degree $>1$ (case i) c)), 
           then $f_1$ does not have compact support, and 

\item [B)] if $\hat\a$ is entire and $\F^{-1}\{\hat\a\}$ has not compact support (case ii)), 
           then $\supp(\A f_1)\subseteq B_1(\vO)$ cannot hold.
\end{itemize}

\subsubsection*{A) $\hat\a$ is a polynomial in $\i\,z$ of degree $>1$}

In the following $\hat\a(z)$ is a polynomial of degree $n>1$ (cf.~(\ref{polynomial})). 
Below we show that the leading term of the polynomial $\hat\a(z)$ is
the most important part for the problem in question and thus we first investigate the 
two cases of monome with even and odd exponent, respectively.

\begin{lemm}\label{lemm:z2n}
Let $N\in\N$, $b\in\R$, $n\in\N^N$ and $z\in\C^N$. If
$\hat\a(z) := b\,z^{2\,n}$, then $f_1$ cannot have compact support.
\end{lemm}

\begin{proof}
We assume that $f_1$ has support in 
\begin{equation}\label{K}
   K:=[-L,L]^N  \qquad\mbox{for some}  \qquad  L>0
\end{equation}
and proof a contradiction. 
Let $z_j=R_j\,e^{\i\,\varphi_j}$ with $R_j>0$ and $\varphi_j\geq 0$ for $1\leq j\leq N$.
For convenience we set
\begin{equation}\label{defRvarphi}
       R:=R_1\,\ldots\,R_N  \qquad \mbox{and}\qquad \varphi:=\sum_{j=1}^N \varphi_j\,.
\end{equation}
From $\hat\a(z) = b\,z^{2\,n}$, it follows 
\begin{equation*}
   \log |f_1(z)|
      = -b\,\mbox{Re}(z^{2\,n})
      = -b\,\mbox{Re}(R^{2\,n}\,e^{\i\,2\,n\,\varphi})
      = -b\,R^{2\,n}\,\cos(2\,n\,\varphi)\,.
\end{equation*}
We have two cases $b>0$ and $b<0$. \\
\begin{itemize}
\item [i)] Let $b>0$. Let $(z^{(m)})_{m\in\N}=((z^{(m)}_1,\,\ldots\,,z^{(m)}_N))_{m\in\N}$ 
      be the sequence such that for $j\in\{1,\,\ldots,\,N\}$
      $$
       (z^{(m)}_j)_{m\in\N}
         = (m\,e^{\i\,\varphi_m})_{m\in\N} \qquad\mbox{with}\qquad
         \varphi_m:= \frac{m\,\pi}{2\,n\,(m+\frac{1}{2})}\,.
      $$
      Then 
      \begin{equation*}
      \log |f_1(z^{(m)})| = -b\,m^{2\,n}\, \cos\left(\frac{m\,\pi}{m+\frac{1}{2}}\right) >0\,
      \end{equation*}
      and (cf. Remark~\ref{rema:PWS} a))
      $$
           H_K(\mbox{Im}(z^{(m)}))= L\,N\,m\,|\sin(\varphi_m)|
           \qquad\mbox{with}\qquad
           \lim_{m\to\infty} \varphi_m=\frac{\pi}{2\,n}\,,
      $$
      where $\sin(\varphi_m)\not=0$.
      Therefore we have for sufficiently large $m$
      \begin{equation}\label{estlogf01}
        \log |f_1(z^{(m)})| \approx C_1\, 
         \left( \frac{H_K(\mbox{Im}(z^{(m)}))}{C_2} \right)^{2\,n}\,,
     \end{equation}
     where $C_1,\,C_2>0$ are constants. We see that condition P2) in Theorem~\ref{th:PWS1} 
     is not satisfied for $K$ defined as in~(\ref{K}) and consequently, $f_1$ does not 
     have compact support for $b>0$.

\item [ii)] Let $b<0$.  For the previous  sequence $(z^{(m)})_{m\in\N}$ in i) with
$\varphi_m:=\varphi:= \frac{\pi}{8\,n}$, we get
\begin{equation*}
      \log |f_1(z^{(m)})| = |b|\,m^{2\,n}\, \cos\left(\frac{\pi}{4}\right) >0 \,
\end{equation*}
which leads again to the estimation~(\ref{estlogf01}). Hence condition P2) in Theorem~\ref{th:PWS1}
is not satisfied and consequently $f_1$ cannot have compact support. This concludes the proof.
\end{itemize}
\end{proof}

\begin{lemm}\label{lemm:z2n-1}
Let $N\in\N$, $b\in\R$, $n\in\N^N$ and $z\in\C^N$. If $\hat\a(z) := \i\,b\,z^{2\,n+1}$, then $f_1$
cannot have compact support.
\end{lemm}

\begin{proof}
We assume that $f_1$ has support in $K$ defined as in~(\ref{K}) and proof a 
contradiction. 
Let $z_j=R_j\,e^{\i\,\varphi_j}$ with $R_j>0$ and $\varphi_j\geq 0$, and $R$, $\varphi$ 
be defined as in~(\ref{defRvarphi}). 
Similarly as in the proof of the previous lemma, it follows that
$$
  \log |f_1(z)| = b\,R^{2\,n+1}\,\sin((2\,n+1)\,\varphi)\,.
$$
Again, we have two cases $b>0$ and $b<0$.
\begin{itemize}
\item [i)] Let $b>0$. Let $(z^{(m)})_{m\in\N}=((z^{(m)}_1,\,\ldots\,,z^{(m)}_N))_{m\in\N}$ 
           be a sequence such that for $j\in\{1,\,\ldots,\,N\}$
           $$
             (z^{(m)}_j)_{m\in\N}
                = (m\,e^{\i\,\varphi_m})_m \qquad\mbox{with}\qquad
               \varphi_m:= \varphi:= \frac{\pi}{2\,(2\,n+1)}\,.
           $$
          Then 
          $$
             \log |f_1(z^{(m)})| = b\,m^{2\,n+1} > 0\,
          $$
          and  (cf. Remark~\ref{rema:PWS} a))
          $$
               H_K(\mbox{Im}(z^{(m)}))= L\,N\,m\,|\sin(\varphi)|\,,
          $$
          where $\sin(\varphi)\not=0$. From this we obtain 
          for sufficiently large $m$ 
          \begin{equation}\label{estlogf02}
             \log |f_1(z^{(m)})| \approx C_1\, 
              \left( \frac{H_K(\mbox{Im}(z^{(m)}))}{C_2} \right)^{2\,n+1}\,,
         \end{equation}
         where $C_1,\,C_2>0$ are constants. Hence property P2) in Theorem~\ref{th:PWS1} 
         is not satisfied for $K$ defined as in~(\ref{K}).

\item [ii)] Let $b<0$. For the previous  sequence $(z^{(m)})_{m\in\N}$ in i) with
            $\varphi:= \frac{3\,\pi}{2\,(2\,n+1)}$, we get the estimation~(\ref{estlogf02}) 
            for sufficiently large $m$, since 
            $$
               b\,\sin((2\,n+1)\,\varphi)>0\,,
            $$ 
            and hence property P2) in Theorem~\ref{th:PWS1} is 
            not satisfied. 
\end{itemize}
Thus, for both cases, condition P2) in Theorem~\ref{th:PWS1} is not satisfied and
therefore $f_1$ cannot have compact support. As was to be shown.
\end{proof}

With the two previous lemmata we now can prove our claim.

\begin{theo}\label{caseP}
If $P(z)$ ($z\in\C^N$) is a polynomial with degree higher than $1$ such that 
$f_1=\F^{-1}\{\exp(-P)\}\in \S(\R^N)$, then $f_1$ cannot have compact support.
\end{theo}

\begin{proof}
Let  $Q_n(z)$ denote the leading term of $-P$ and $g_n:= \exp\{Q_n\}$. Then the degree $n$
of the leading term is larger than $1$ and according to Lemma~\ref{lemm:z2n} if $n$ is even
or Lemma~\ref{lemm:z2n-1} if $n$ is odd, there exists a sequence $(z^{(m)})_{m\in\N}$ such that
$g_n(z^{(m)})$ satisfies the estimation~(\ref{estlogf01}) or~(\ref{estlogf02}) with $f$ 
replaced by $g_n$.
This together with\footnote{This line of argumentation cannot be used for case B) below, where 
$P\equiv\hat\a$ is a power series.}
\begin{equation*}
\begin{aligned}
       \mbox{Re}((-P-Q_n)(z^{(m)})) <<  \mbox{Re}(Q_n(z^{(m)})) \,
\qquad\mbox{for sufficiently large $m$,}
\end{aligned}
\end{equation*}
implies that $f_1(z^{(m)})$ satisfies the estimation~(\ref{estlogf01}) or~(\ref{estlogf02}). 
And thus condition P2) in Theorem~\ref{th:PWS1} is not satisfied, i.e.
$f_1$ cannot have compact support. This proves the Theorem.
$\Box$
\end{proof}

\subsubsection*{B) $\hat\a$ is entire and $\F^{-1}\{\hat\a\}$ has not compact support}

For the proof we need an estimation of the minimum modulus of
entire functions. Theorem 3.7.4 in~\cite{Boa73} fulfills our requirements.
We state it for the case of functions of exponential type $\tau$.
We recall that $F:\C\to\C$ ($N=1$) is of exponential type $\tau$ if
\begin{equation*}
\begin{aligned}
      \limsup_{r\to\infty} \frac{\log M(r)}{r} = \tau\,,
\end{aligned}
\end{equation*}
where 
$M(r):=\max_{|z|=r}|F(z)|$ (cf. Definition 2.1.3 in~\cite{Boa73}).

\begin{theo}\label{MinModulusTh}
If the entire function $F:\C\to\C$ is of exponential type $\tau$, then for every
$\epsilon,\,\eta>0$, $\xi>1$ and sufficiently large $R>0$ we have
\begin{equation*}
    \log |F(z)| > -A\,(\tau+\epsilon)\,\xi\,R
\qquad \mbox{for} \qquad |z|<R\,
\end{equation*}
except in a set of circles the sum of whose radii is at most $2\,\eta\,\xi\,R$.
Here $A>0$ depends only on $\eta$ and $\xi$.
\end{theo}

Let $\hat a(z)$ be entire and $F(z):=\exp\{-\hat a(z)\}$ be a function of exponential 
type $\tau=1$. Moreover, let $w(z_1)$ for $z_1\in\C$ be defined as in~(\ref{defwz}). Then 
Theorem~\ref{MinModulusTh} implies for the setting $\epsilon\in (0,1)$, $\eta=1/8$, $\xi=2$,
$\tau=1$ that
$$
     \mbox{Re}(\hat a(w(z_1)))  <  4\,A\,R
      \qquad \mbox{for} \qquad |z_1|=3\,R/4 > 2\,\eta\,\xi\,R\,.
$$
We emphasize that the constant $A$ does not depend on $R$ and therefore $R$ can be replaced
by any larger number $\tilde R$. For the following proof we use the reformulation
\begin{equation}\label{MMProp}
\begin{aligned}
    |\exp\{-\hat a(w(z_1))\}|
      &=  \exp\{-\mbox{Re}(\hat a(w(z_1)))\}  \\
      &>  \exp\{-B\,r\}
      \qquad\qquad\qquad \mbox{for $|z_1|=r$}\,
\end{aligned}
\end{equation}
for sufficiently large $r>0$.  Here $B:=16\,A/3>0$ and $r\equiv 3\,R/4$.

\begin{theo}\label{theo:caseSeries}
Let $\hat\a$ and $\A$ be defined as in Definition~\ref{defaA} and $G$ be defined as
in~(\ref{defG}). If $\hat\a$ is entire, $\F^{-1}\{\hat\a\}$ has not compact support
and~(\ref{A2}) holds, then $\supp(\A f_1)\subseteq B_1(\vO)$ cannot hold.
\end{theo}

\begin{proof}
We perform a proof by contradiction and assume that $\supp(\A f_1)\subseteq B_1(\vO)$ holds.  
Let $w(z_1)$ for $z_1\in\C$  be defined as in~(\ref{defwz}). Then the assumption together 
with Theorem~\ref{th:PWS1} and Remark~\ref{rema:PWS} b) imply that
\begin{equation}\label{FAfw}
   \F\{\A\,f_1\}\circ w:\C\to\C,\,z_1\mapsto \hat\a(w(z_1))\,\exp\{-\hat\a(w(z_1))\}
\end{equation}
is an exponential function of type $\leq 1$. We prove the contradiction by showing 
that $\F\{\A\,f_1\}\circ w$ cannot be of exponential type $\tau\in (0,\infty)$. 
Since, by assumption, $\hat\a(z)$ is entire, rotational symmetric and $\F^{-1}\{\hat\a\}$ 
does not have compact support, Remark~\ref{rema:PWS} c) and a) imply that for each $C>0$ 
there exists an $m\in\N$ such that 
\begin{equation*}
\begin{aligned}
      \max_{|z_1|=m}  |\hat a(w(z_1))|
            > C\,\max_{|z_1|=m} \left(\exp\{ m\,|\mbox{Im}(z_1)| \} \right)  
            = C\,\exp\{ m^2 \}  \,. 
\end{aligned}
\end{equation*}
Let $M(r)$ denote the maximum modulus function of $\F\{\A f_1\}(w(z_1))$  
on the circle $|z_1|=r$.
From the previous estimation,~(\ref{FAfw}) and~(\ref{MMProp}), we get
\begin{equation*}
\begin{aligned}
   M(m) 
       &= \max_{|z_1|=m} \left| \hat\a(w(z_1))\,\exp\{-\hat\a(w(z_1))\} \right|\\
       &> C\,\exp\{ m\,\left(m - B\right)\}\,
\end{aligned}
\end{equation*}
for sufficiently large $m$. 
Without loss of generality we assume $C=1$. Then we obtain
\begin{equation*}
   \limsup_{m\to\infty} \frac{\log M(m)}{m}
      \geq \limsup_{m\to\infty} \,(m -B)
      =   \infty\,
\end{equation*}
and thus we have proven that $\F\{\A\,f\}(w(z_1))$ cannot be an exponential function of
type $\tau\in (0,\infty)$ which concludes the proof. 
\end{proof}


\begin{thebibliography}{10}

\bibitem{Boa73}
Boas, R. P.
\newblock {\em Entire Functions}.
\newblock Academic Press Inc., New York, 3. Printing, 1973.

\bibitem{DauLio92_5}
Dautray, R. and Lions J.-L.
\newblock {\em Mathematical Analysis and Numerical Methods for Science and Technology. Volume 5.}
\newblock {\em Springer-Verlag}, New York, 1992.

\bibitem{EngRun95}
Engl, H. W. and Rundell, W. (eds.),
\newblock {\em Inverse Problems in Diffusion Processes}.
\newblock SIAM, Philadelphia, 1995.

\bibitem{EngHanNeu96}
Engl, H.W. and Hanke, M. and Neubauer, A.
\newblock {\em Regularization of Inverse Problems}.
\newblock Kluwer Academic Publishers, Dordrecht, 1996.

\bibitem{EVa99}
Evans, L. C.
\newblock {\em Partial Differential Equations}.
\newblock American Mathematical Society, Providence, Rhode Island, 1999.

\bibitem{Fet80}
Fetter, A.L. and Walecka, J.D.
\newblock {\em Theoretical Mechanics of Particles and Continua}. 
\newblock McGraw-Hill, New York, 1980. 

\bibitem{GasWit99}
Gasquet, C. and Witomski. P.
\newblock {\em Fourier Analysis and Applications}.
\newblock Springer Verlag, New York,  1999.

\bibitem{GuiMorRya04}
Guichard, F. and Morel, J.-M. and Ryan, R.
\newblock {\em Contrast invariant image analysis and PDE's"}
\newblock Lecture notes, see http://mw.cmla.ens-cachan.fr/~morel/, 2004.

\bibitem{Harris79}
Harris, C. J.
\newblock {\em Mathematical Modelling of Turbulent Diffusion in the Environment}. 
\newblock Academic Press, New York, 1979. 

\bibitem{Hoe03}
H\"ormander, L.
\newblock {\em The Analysis of Linear Partial Differential Operators I}.
\newblock Springer Verlag, New York, 2nd edition, 2003.

\bibitem{Isa98}
Isakov, V. 
\newblock {\em Inverse Problems for Partial Differential Equations}.
\newblock Springer Verlag, New York, 1998.

\bibitem{KiSrTr06}
Kilbas, A. A. and Srivastava, H. M. and Trujillo, J. J. 
\newblock {\em Theory and Applications of Fractional Differential Equations}.
\newblock Elsevier, New York, 2006.

\bibitem{Kir96}
Kirsch, A. 
\newblock {\em An Introduction to the Mathematical Theory of Inverse Problems}.
\newblock Springer Verlag, New York, 1996.

\bibitem{KitKro93}
Kittel, Ch. and {Kr\"omer}, H.
\newblock {\em Einf\"uhrung in die Festk\"orperphysik}.
\newblock  Oldenbourg Verlag, M\"unchen,  9. Auflage, 1991.

\bibitem{KowSchBon10}
Kowar, R. and Scherzer, O. and Bonnefond, X.
\newblock {Causality analysis of frequency-dependent wave attenuation.}
\newblock {Math. Meth. Appl. Sci.} 2010, DOI: 10.1002/mma.1344

\bibitem{Kow10}
Kowar, R.
\newblock {Integral equation models for thermoacoustic imaging of acoustic
                 dissipative tissue.}
\newblock {Inverse Problems 26 (2010), 095005 (18pp)},
DOI: 10.1088/0266-5611/26/9/095005

\bibitem{KowSch10}
Kowar, R. and Scherzer, O.
\newblock {Photoacoustic Imaging Taking into Account Attenuation.}
(arXiv:1009.4350)
\newblock {to appear in: \newblock{Mathematical Modeling
                        in Biomedical Imaging II: Optical, Ultrasound, and
                        Opto-Acoustic Tomographies,  Lecture Notes in
                        Mathematics: Mathematical Biosciences Subseries}},


\bibitem{NacSmiWaa90}
A.~I.~Nachman, J.~F.~Smith, III and R.~C.~Waag.
\newblock An equation for acoustic propagation in inhomogeneous media with relaxation losses.
\newblock {\em J. Acoust. Soc. Am.} 88 (3), Sept. 1990.

\bibitem{Rem92a}
Remmert, R.
\newblock {\em Funktionentheorie 1}.
\newblock Springer Verlag, New York, 2.Auflage, 1992.

\bibitem{Rem95b}
Remmert, R.
\newblock {\em Funktionentheorie 2}.
\newblock Springer Verlag, New York, 3.Auflage, 1995.

\bibitem{ScGrGrHaLe09}
Scherzer, O. and Grasmair, M. and Grossauer, H. and Haltmeier, M. and Lenzen, F.
\newblock {\em  Variational Methods in Imaging}.
\newblock Springer-Verlag, New York, 2009.

\bibitem{Soi99}
Soille, P.
\newblock {\em Morphological image analysis. Principles and applications}. 
\newblock Springer-Verlag, Berlin, 1999.

\bibitem{Wei98a}
Weickert, J.
\newblock {\em Anistropic Diffusion in Image Processing}.
\newblock Teubner Stuttgart Verlag, Stuttgart, 1998.

\bibitem{Wei72}
Weinberg, S.
\newblock {\em Gravitation and Cosmology: Principles and Applications of the General 
           Theory of Relativity}.
\newblock John Wiley \& Sons, New York, 1972.

\bibitem{Wei95}
Weinberg, S.
\newblock {\em The Quantum Theory of Field. Vol. I. Foundations}.
\newblock Cambridge University Press, Cambridge, 1995.



\bibitem{Zee03}
Zee, A.
\newblock {\em Quantum Field Theory in a Nutshell}.
\newblock Princeton University Press, New Jersey, 2003.



\end{thebibliography}
\end{document}